\documentclass[11 pt, psamsfonts, leqno]{amsart}
\usepackage{amssymb}
\usepackage{amsmath}	

\usepackage{bm}

 \calclayout
\def\ignore#1{\relax}

\ignore{
\usepackage[scaled=0.92]{helvet}	
\usepackage[subscriptcorrection,slantedGreek,nofontinfo,  mtphbi]{mtpro2}
}

\usepackage[plainpages=false, pdfpagelabels,  colorlinks=true, pdfstartview=FitV, linkcolor=blue, citecolor=blue, urlcolor=blue]{hyperref}

\usepackage{epic}
\usepackage{eepicemu}

\usepackage{graphicx}
\usepackage{epstopdf}
\usepackage{maxidiagram}

\newcommand\inlinegraphic[2][{scale=1.0}]{\begin{array}{c} \includegraphics[#1]{./EPS/#2}\end{array}}

\numberwithin{equation}{section}
\numberwithin{figure}{section}

\def\Z{{\mathbb Z}}

\def\Q{{\mathbb Q}}
\def\C{{\mathbb C}}

\def\eps{\varepsilon}
 
 \def\S{{\mathfrak{S}}}

\def\spp #1{^{(#1)}}

\def\u #1 #2{\mathcal U(#1, #2)}  
\def\uhat #1 #2{\widehat{\mathcal U}(#1, #2)}  

\def\kt #1{{KT_{#1}}}  
\def\akt #1{\widehat{KT}_{#1}}  
\def\abmw #1{\widehat{W}_{#1}}  
\def\bmw #1{W_{#1}}  
\def\ahec #1{\widehat{H}_{#1}}  
\def\hec #1{H_{#1}}
\def\w #1 #2{\bmw {#1}^{(#2)}}  
\def\V #1 #2{V_{#1}^{(#2)}}  
\def\k #1 #2{ KT_{#1}^{(#2)}}
\def\br{D}   
\def\br #1{D_{#1}} 
\def\p #1{\bm {#1}}
\def\pbar #1{\bar{\p #1}}

\def\inv{^{-1}}
\def\la{\lambda}
\def\vtbold{\varthetabold}
\def\labold{\lambdabold}
\def\vt{\vartheta}
\def\cl{{\rm cl}}

\def\rhobold{{\bm \rho}} 
\def\deltabold{{\bm \delta}}

\def\qbold{{\bm q}}
\def\ubold{{\bm u}}
\def\varthetabold{{\bm \vartheta}}
\def\lambdabold{{\bm \lambda}}
\def\mubold{{\bm\mu}}
\def\powerpm{^{\pm 1}} 
\def\abold{\bm a}
\def\A{\mathbb A}
\def\B{\mathbb B}

\def\unitvector #1{\bm{ \hat {#1}}}

\def\End{{\rm End}}

\def\Res{{\rm Res}}
\def\id{{\rm id}}

\def\hods{\unskip\kern.55em\ignorespaces}

\theoremstyle{plain}
\newtheorem{theorem}{Theorem}[section]

\theoremstyle{plain}
\newtheorem{proposition}[theorem]{Proposition}

\theoremstyle{plain}
\newtheorem{corollary}[theorem]{Corollary}

\theoremstyle{plain}
\newtheorem{lemma}[theorem]{Lemma}

\theoremstyle{definition}
\newtheorem{definition}[theorem]{Definition}
\theoremstyle{definition}

\theoremstyle{definition}
\newtheorem{remark}[theorem]{Remark}

\theoremstyle{remark}

\title[Cyclotomic BMW algebras, II]{Cyclotomic Birman--Wenzl--Murakami algebras,  II:\\ Admissibility Relations and Freeness }

\author{Frederick M. Goodman}
\address{ Department of Mathematics\\ University of Iowa\\ Iowa
City, Iowa}
\email{ goodman@math.uiowa.edu} 
\author{Holly Hauschild Mosley}
\address{Department of Mathematics\\ Grinnell College \\
Grinnell, Iowa}
\email{HAUSCHIL@GRINNELL.EDU}


\begin{document}
 \baselineskip=16pt 
 \begin{abstract}  The cyclotomic Birman-Wenzl-Murakami algebras are quotients
of the affine BMW algebras in which the affine generator satisfies a polynomial relation.  
We study admissibility conditions on the ground ring for these algebras, and show
that the algebras defined over an admissible integral ground ring $S$ are free $S$--modules and isomorphic to cyclotomic  Kauffman tangle algebras.
We also determine the representation theory in the generic semisimple case, obtain a recursive formula for the weights of the Markov trace,  and give a sufficient condition for semisimpliity.
\end{abstract}

 \maketitle

March, 2008.

\setcounter{tocdepth}{1}
\tableofcontents

\section{Introduction} 

This paper and the companion paper ~\cite{GH2} continue the study of affine and cyclotomic Birman--Wenzl-- Murakami  (BMW) algebras, which we began in  ~\cite{GH1}.  

\subsection{Background} 
The origin of the BMW algebras was in knot theory.  Kauffman defined ~\cite{Kauffman} an invariant of regular isotopy for  links in $S^3$, determined by certain skein relations.
 Birman and Wenzl ~\cite{Birman-Wenzl} and independently Murakami ~\cite{Murakami-BMW}  then defined  a family of quotients of the braid group algebras, and showed that Kauffman's invariant could be recovered from a trace on these algebras.  These (BMW) algebras were  defined by generators and relations, but were implicitly modeled on certain algebras of tangles whose definition  was subsequently made explicit by Morton and Traczyk ~\cite{Morton-Traczyk}, as follows:
  Let $S$ be a commutative unital ring with invertible elements
$\rho$, $q$, and $\delta_0$ satisfying $\rho\inv - \rho = (q\inv -q) (\delta_0 - 1)$.  The {\em Kauffman tangle algebra}  $\kt{n, S}$  is the $S$--algebra of framed $(n, n)$--tangles in the disc cross the interval,  modulo Kauffman skein relations:
\begin{enumerate}
\item Crossing relation:
$
\quad \inlinegraphic[scale=.6]{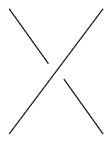} - \inlinegraphic[scale=.3]{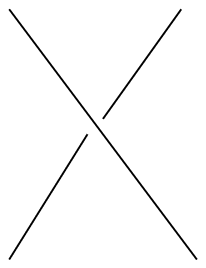} 
\quad = 
\quad
(q\inv - q)\,\left( \inlinegraphic[scale=1]{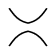} - 
\inlinegraphic[scale=1]{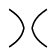}\right).
$
\item Untwisting relation:
$\quad 
\inlinegraphic{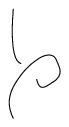} \quad = \quad \rho \quad
\inlinegraphic{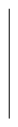} \quad\ \text{and} \quad\ 
\inlinegraphic{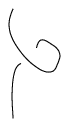} \quad = \quad \rho\inv \quad
\inlinegraphic{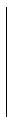}. 
$
\item  Free loop relation:  $T\, \cup \, \bigcirc = \delta_0 \, T. $
\end{enumerate}
Morton and Traczyk ~\cite{Morton-Traczyk}   showed that the $n$--strand algebra $\kt{n, S}$ is free of rank $(2n-1)!!$ as a module over $S$, and
Morton and Wassermann ~\cite{Morton-Wassermann} proved that the BMW algebras and the Kauffman tangle algebras are isomorphic.
 
It is natural to ``affinize" the BMW algebras to obtain BMW analogues of the affine Hecke   algebras of type $A$, see ~\cite{ariki-book}.  The affine Hecke algebra  can be realized geometrically as the algebra of braids in the annulus cross the interval,  modulo Hecke skein relations; this suggests defining the affine Kauffman tangle algebra $\akt{n, S}$ as the algebra of framed $(n, n)$--tangles in the annulus cross the interval, modulo Kauffman skein relations.  Turaev  ~\cite{Turaev-Kauffman-skein} showed that the resulting algebra of  $(0,0)$--tangles is a  (commutative) polynomial algebra in infinitely many variables, so it  makes sense  to absorb this polynomial algebra into the ground ring.  (The ground ring gains infinitely many parameters $\delta_j$ ($j \ge 1$)  corresponding to the generators of the polynomial algebra.)
On the other hand,  H\"aring--Oldenburg ~\cite{H-O2} defined an affine version of the BMW algebras by generators and relations.  In ~\cite{GH1}, we showed that 
H\"aring--Oldenburg's affine BMW algebras are isomorphic to the affine Kauffman tangle algebras,  and we showed that these algebras are free modules over their ground ring, with a basis reminiscent of a well--known basis of affine Hecke algebras.

The affine BMW algebras arise naturally in several different contexts:  knot theory in the solid torus ~\cite{Turaev-Kauffman-skein, Lambr-thesis, Lambr-solid-torus},  representations of Artin braid group of type $B$ by $R$ matrices of quantum groups ~\cite{RO},  and representations of ordinary BMW algebras.   See ~\cite{GH2},  Section 1.1
 for more detail.

\subsection{Cyclotomic algebras}
 In this paper and the companion paper ~\cite{GH2} we consider cyclotomic BMW algebras,  which are the BMW analogues of cyclotomic Hecke algebras ~\cite{ariki-book}.  The affine BMW algebras have a distinguished generator $x_1$, which, in the geometric (Kauffman tangle) picture is represented by a braid with one strand wrapping around the hole in the annulus cross interval.  The {\em cyclotomic BMW algebra} $\bmw{n, S, r}$ is defined to be the quotient of the affine BMW algebra $\abmw{n, S}$ in which the generator $x_1$ satisfies a monic polynomial equation
\begin{equation} \label{cyclotomic polynomial relations 0}
x_1^r + \sum_{k=0}^{r-1} a_k x_1^k = 0.
\end{equation}
with coefficients in $S$.\footnote{Actually, we will assume that  the polynomial splits in $S$.}  The cyclotomic BMW algebras were also introduced by H\"aring-Oldenburg in ~\cite{H-O2}.

In the geometric (Kauffman tangle) picture, it is more natural to convert this relation into a local skein relation:
\begin{equation} \label{cyclotomic skein relations 0}
 T_r + \sum_{k=0}^{r-1} a_k  T_k = 0,
\end{equation}
whenever $T_0, T_1, \dots, T_r$  are affine tangle diagrams that are identical in the exterior of some disc $E$ and    $T_k \cap E$  consists of one strand wrapping $k$ times around the hole in the annulus cross interval; i.e.  $T_k \cap E$ ``equals" $x_1^k$.    The {\em cyclotomic Kauffman tangle algebra} $\kt{n, S, r}$  is defined to be the quotient of the affine Kauffman tangle algebra
$\akt{n, S}$ by the cyclotomic skein relation.  See ~\cite{GH2},  Section 1.2,
for an alternative description of the cyclotomic Kauffman tangle algebras in terms of the affine Kauffman tangle category.

A priori, the ideal  in $\akt{n, S} \cong \abmw{n,S}$ generated by the cyclotomic skein relation (\ref{cyclotomic skein relations 0}) is larger than the ideal generated by the polynomial relation (\ref{cyclotomic polynomial relations 0}), so we have a surjective, but not evidently injective homomorphism $$\varphi : \bmw{n, S, r} \to \kt{n, S, r}.$$

\subsection{Admissibility}  The cyclotomic BMW algebras and Kauffman tangle algebras can be defined over an arbitrary commutative unital ring $S$ with parameters $\rho$, $q$,  $\delta_j$    ($j \ge 0$), and
$u_1, \dots, u_r$ (the eigenvalues of $y_1$), assuming   that $\rho$, $q$,  $\delta_0$  and  $u_1, \dots, u_r$   are invertible, 
and $\rho\inv - \rho=   (q\inv -q) (\delta_0 - 1)$.    However, unless the parameters
satisfy additional relations,  the identity element $\bm 1$ of the cyclotomic Kauffman tangle algebras
will be a torsion element over $S$;  if $S$ is a field (and the additional relations do not hold)  then $\bm 1 = 0$,  so $\kt {n, S, r}$ becomes trivial.    The additional conditions are called ``weak admissibility;"  see Section \ref{subsection: weak admissibility}.   Weak admissibility is thus a minimal condition for the non--triviality of the cyclotomic algebras.

In order to obtain substantial results about the cyclotomic BMW and Kauffman tangle algebras, it seems necessary to impose a  condition on the ground ring $S$ that is stronger than weak admissibility.
An appropriate condition was introduced by Wilcox and Yu in ~\cite{Wilcox-Yu}.  
Their condition has a simple formulation in terms of the 
 representation theory of the 2--strand algebra $\bmw{2, S, r}$,  and also translates    into explicit relations on the parameters.  The condition of Wilcox and Yu,   called ``admissibility,"   is the subject of 
Section \ref{section:  representations and admissibility}  of this paper.

\subsection{Results}  The main results of this paper and ~\cite{GH2} is that if the ground ring $S$ is an integral domain and admissible in the sense of Wilcox and Yu,  then $\bmw{n, S, r} \cong \kt{n, S, r}$,  and, moreover, these algebras are free $S$--modules of rank $r^n (2n-1)!!$.  The proof of these results has a topological component,  which is given in ~\cite{GH2},  and an algebraic component, which is given in this paper.

The topological arguments in ~\cite{GH2} have two essential consequences:  first they show that  the cyclotomic BMW algebra $\bmw{n,S,r}$  has a spanning set $\A'_r$  of cardinality $r^n (2n-1)!!$.
Second, they show that the  $(0,0)$--tangle algebra $\kt{0, S, r}$  is free of rank 1, assuming that 
$S$ is weakly admissible.  This implies the existence of the ``Markov trace"  on the cyclotomic BMW and Kauffman tangle algebras (and, hence, Kauffman type invariants for links in the solid torus).

To achieve the main results, it remains to show that the spanning set $\varphi(\A'_r)$
of the cyclotomic Kauffman tangle algebra $\kt{n, S, r}$ is linearly independent, when the ground ring $S$ is an admissible integral domain.

To do this we first show that there exists a universal admissible integral domain $\overline R_+$ with the property that every admissible integral domain is a quotient of $\overline R_+$.  Moreover,  field of fractions
$F$ of $\overline R_+$ is the field of rational functions  $\Q(\qbold, \ubold_1, \dots, \ubold_r)$.

\renewcommand\mod{\rm{\ mod\ }}

We then analyze the representation theory of the cyclotomic BMW algebras   defined over  $F$, by adapting the inductive method of Wenzl from ~\cite{Birman-Wenzl, Wenzl-Brauer, Wenzl-BCD}.   This analysis shows that $\kt{n, F, r} \cong \bmw{n, F, r}$, and that these algebras are split  simisimple of dimension  $r^n (2n-1)!!$.   The irreducible representations of $\bmw{n, F, r}$ are indexed by $r$-tuples of  Young diagrams of total size $\le n$ and congruent to $n \mod 2$.
\ignore{A novel feature of our argument is that we incorporate into the inductive proof both the isomorphism of the cyclotomic BMW and Kauffman tangle algebras and recursive formulas for the weights of the Markov trace.}   

A relatively simple argument   then shows that for any admissible integral domain $S$,  we have
$\kt{n, S, r} \cong \bmw{n, S, r}$, and these algebras are free of rank $r^n (2n-1)!!$ over $S$. 

As a biproduct of our argument, we obtain a sufficient condition for semisimplicity of the cyclotomic BMW algebras defined over an admissible field.  

The elements
$$
y_j =  g_{j-1}  \cdots     g_1  y_1 g_1 \cdots g_{j-1},
$$
are important for  the representation theory of the cyclotomic BMW algebras.  They are BMW analogues of Jucys-Murphy elements; 
see, for example,  \cite{leduc-ram, mathas-JM}.  Indeed, in our treatment of the representation theory over the generic field $F$,  the action of the Jucys-Murphy elements plays an essential role.

\subsection{Related work, and acknowledgments} Wilcox and Yu  have been studying the same material independently and have obtained  similar (and, in some instances, slightly stronger)  results ~\cite{Wilcox-Yu, Yu-thesis, Wilcox-Yu2}.   We are indebted to 
 Wilcox and Yu for pointing out an error in a previous preprint version of our work,  which required us to substantially rework the algebraic (linear independence) component of the arguments.  In fact, we had to adopt a completely different strategy for proving linear independence.  Meanwhile,  Wilcox and Yu ~\cite{Yu-thesis, Wilcox-Yu2} were able to make our original strategy work, using  a refined analysis of their admissibility condition. 

We would also like to mention here recent work of 
Ariki, Mathas and Rui on the 
``cyclotomic Nazarov--Wenzl algebras" ~\cite{ariki-mathas-rui},  which are cyclotomic quotients of the degenerate affine BMW algebras introduced by Nazarov ~\cite{Nazarov}.
The idea of deriving an admissibility condition on the ground ring from an assumption on the representation theory of the 2--strand algebra originates in  ~\cite{ariki-mathas-rui}, and is essential to our work and the work of Wilcox and Yu.

We rely heavily on ideas from Beliakova and Blanchet   ~\cite{blanchet-beliakova}  (and indirectly on Nazarov ~\cite{Nazarov}) for our computation of the weights of the Markov trace and of the action of conjugates of the affine generator $y_1$  on the  basis of up--down tableaux.

It is shown in ~\cite{goodman-2008}  and ~\cite{Yu-thesis}  that cyclotomic BMW algebras defined over integral admissible ground rings are cellular
in the sense of Graham and Lehrer ~\cite{Graham-Lehrer-cellular}.  The proof in ~\cite{goodman-2008} depends on the isomorphism of cyclotomic BMW algebras and cyclotomic Kaufmann tangle algebras established here.

Finally,  a recent preprint of Rui and Xu  ~\cite{rui-2008} concerns the representation theory of cyclotomic BMW algebras 
over a $u$--admissible field.

\section{Preliminaries} 

 \subsection{The affine BMW algebra} 
  
We recall the definition of  the affine Birman--Wenzl--Murakami  (BMW) algebra.  Let
 $S$ be a commutative unital ring with
  elements 
$\rho$, $q$, and  $\delta_j$  ($j \ge 0$),   such that $\rho$, $q$, and  $\delta_0$ are invertible, and  the relation
$
\rho\inv - \rho=   (q\inv - q) (\delta_0 - 1)
$
holds.

\begin{definition} \label{definition affine BMW}
 The {\em affine
Birman--Wenzl--Murakami} algebra
$\abmw  {n, S}$ is the
$S$ algebra with generators $y_1^{\pm 1}$, $g_i^{\pm 1}$  and
$e_i$ ($1 \le i \le n-1$) and relations:
\begin{enumerate}
\item (Inverses)\hods $g_i g_i\inv = g_i\inv g_i = 1$ and 
$y_1 y_1\inv = y_1\inv y_1= 1$.
\item (Idempotent relation)\hods $e_i^2 = \delta_0 e_i$.
\item (Affine braid relations) 
\begin{enumerate}
\item[\rm(a)] $g_i g_{i+1} g_i = g_{i+1} g_ig_{i+1}$ and 
$g_i g_j = g_j g_i$ if $|i-j|  \ge 2$.
\item[\rm(b)] $y_1 g_1 y_1 g_1 = g_1 y_1 g_1 y_1$ and $y_1 g_j =
g_j y_1 $ if $j \ge 2$.
\end{enumerate}
\item[\rm(4)] (Commutation relations) 
\begin{enumerate}
\item[\rm(a)] $g_i e_j = e_j g_i$  and
$e_i e_j = e_j e_i$  if $|i-
j|
\ge 2$. 
\item[\rm(b)] $y_1 e_j = e_j y_1$ if $j \ge 2$.
\end{enumerate}
\item[\rm(5)] (Affine tangle relations)\vadjust{\vskip-2pt\vskip0pt}
\begin{enumerate}
\item[\rm(a)] $e_i e_{i\pm 1} e_i = e_i$,
\item[\rm(b)] $g_i g_{i\pm 1} e_i = e_{i\pm 1} e_i$ and
$ e_i  g_{i\pm 1} g_i=   e_ie_{i\pm 1}$.
\item[\rm(c)\hskip1.2pt] For $j \ge 1$, $e_1 y_1^{ j} e_1 = \delta_j e_1$. 
\vadjust{\vskip-
2pt\vskip0pt}
\end{enumerate}
\item[\rm(6)] (Kauffman skein relation)\hods  $g_i - g_i\inv = (q - q\inv)(1- e_i)$.
\item[\rm(7)] (Untwisting relations)\hods $g_i e_i = e_i g_i = \rho \inv e_i$
 and $e_i g_{i \pm 1} e_i = \rho  e_i$.
\item[\rm(8)] (Unwrapping relation)\hods $e_1 y_1 g_1 y_1 = \rho e_1 = y_1 
g_1 y_1 e_1$.
\end{enumerate}
\end{definition}

\begin{remark}  The presentation differs  slightlyÊ\ from the one we used in ~\cite{GH1}.    There we used the generator $x_1 = \rho \inv y_1$,   and
 parameters $\vartheta_j =  \rho^{-j} \delta_j$    (so that $e_1 x_1^j e_1 = \vartheta_j e_1$).   We also used the parameter $z$  in place of $(q\inv -q)$.
   \end{remark}

In $\abmw{n, S}$,  define  for $1 \le j \le n$,
$$
y_j =  g_{j-1}  \cdots     g_1  y_1 g_1 \cdots g_{j-1}
$$
and
$$
y_j' =  g_{j-1}  \cdots     g_1  y_1 g_1\inv \cdots g_{j-1}\inv.
$$

We record some elementary results from ~\cite{GH1}:

\begin{lemma}\label{proposition-xr relations1}
The following relations hold in $\abmw{n, S}$.   For $1 \le k, j \le n$,
\begin{enumerate}
\item For $j \not\in  \{k, k-1\}$,  $g_j  y_k =  y_k g_j$. 
\item  For  $j \not\in  \{k, k-1\}$,   
$e_j  y_k =  y_k e_j$.
\item   $y_j y_k = y_k y_j$.
\item  For $k < n$, $g_k y_k = y_{k+1} g_k \inv$ and $g_k \inv y_{k+1} = y_k 
g_k$.
\item For $k < n$, 
\begin{eqnarray*}
 g_k y_{k+1} &=& y_k g_k  + (q - q\inv) y_{k+1} - (q - q\inv) \rho e_k y_{k+1},\\
g_k \inv y_k &=& y_{k+1} g_k\inv  - (q - q\inv) y_k  + (q - q\inv) e_k y_k.
\end{eqnarray*}
\item  For $k < n$, $e_k  y_k =   e_k y_{k+1}\inv$ and 
 $e_{k} y_k \inv =   e_{k} y_{k+1} $.
\item  $g_{k}y_{k}g_{k} y_{k} = y_{k}g_{k} y_{k} g_{k}$.
\item  $e_{k} y_{k} g_{k} y_{k} = \rho\, e_{k}$.
\item For $s \ge 1$,  $e_{k} y_{k}^{s} g_{k} y_{k} =   e_{k} y_{k}^{s-1} g_{k}\inv$, and 
 $e_{k} y_{k}^{-s} g_{k}\inv y_{k}\inv =  e_{k} y_{k}^{-s+1} g_{k}$.

\end{enumerate}
\end{lemma}

\begin{remark}  \label{remark: involution on affine BMW}
There is a unique antiautomorphism of the affine BMW algebra leaving each of the generators $y_1,  e_i, g_i$ invariant; the antiautomorphism simply reverses the order of a word in the generators.    This follows from the symmetry of the relations.
Using this, one obtains versions of points (6), (8) and (9) of the Lemma with $e_k$ written on the right,  for example:
$
 y_{k}   g_{k} y_{k}^{s} e_{k}  =   g_{k}\inv y_{k}^{s-1}  e_{k}.
$
\end{remark}

\begin{lemma}  \label{lemma - recursion for f_r}
 For $j \ge 1$,  there exist elements $\delta_{-j} \in \Z[\rho^{\pm 1}, q^{\pm 1}, \delta_0, \dots, \delta_j]$ such that $e_1 y_1^{-j} e_1 = \delta_{-j} e_1$.  Moreover,  the elements $\delta_{-j}$  are determined by the recursion relation:
\begin{equation} \label{equation:  delta(-j)  recursive relations}
\begin{aligned}
\delta_{-1} &= \rho^{-2} \delta_1 \\
\delta_{-j}  & =  \rho^{-2}  \delta_j +  (q\inv - q) \rho\inv  \sum_{k = 1}^{j-1} (\delta_k \delta_{k-j} - \delta_{2k -j}   ) \quad  (j \ge 2).
\end{aligned}
\end{equation}

\begin{proof}  Follows from ~\cite{GH1},   Corollary 3.13, and
~\cite{GH2}, Lemma 2.6.
\end{proof}

\end{lemma}

\subsection{The affine Kauffman tangle algebra.} 
 Let
 $S$ be a commutative unital ring with
  elements 
$\rho$, $q$, and  $\delta_j$  ($j \ge 0$),   such that $\rho$, $q$, and  $\delta_0$ are invertible, and  the relation
$
\rho\inv - \rho=   (q\inv - q) (\delta_0 - 1)
$
holds.
The {\em affine Kauffman tangle algebra}  $\akt{n, S}$    is the $S$--algebra
of framed $(n, n)$ tangles in  $A \times I$, where $A$ is the annulus and $I$ the interval,  modulo the Kauffman skein relations.     The reader is referred to 
 ~\cite{GH1}  for a detailed definition.
 
 The affine Kauffman tangle algebra  $\akt{n, S}$ can be described in terms of 
 {\em affine tangle diagrams}, which  are quasi-planar diagrams representing framed tangles.  Here is a picture of a typical affine tangle diagram.  The heavy vertical line represents the hole in  $A \times I$; we refer to it as the flagpole.
$$
\inlinegraphic[scale=1.5]{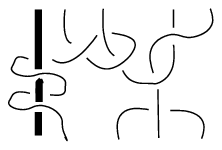}
$$
Multiplication is by stacking of such diagrams; our convention is that the product $a b$ is given by stacking $b$ over $a$.    The affine Kauffman tangle algebra is generated by the following affine tangle diagrams:
$$
X_1 = \inlinegraphic[scale= .7]{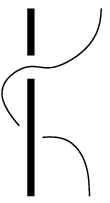}
\qquad
G_i =  \inlinegraphic[scale=.6]{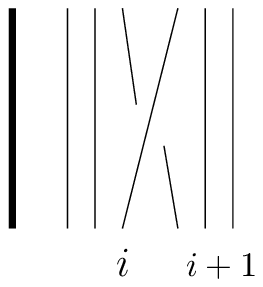}\qquad
E_i =  \inlinegraphic[scale= .7]{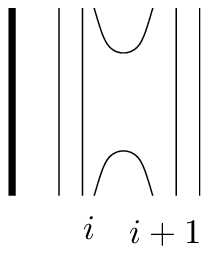} 
$$

\begin{theorem}[\cite{GH1}]  The affine BMW algebra $\abmw{n, S}$  is isomorphic to  the affine Kauffman tangle algebra    $\akt{n, S}$ by a map determined by the assignments
$y_1 \mapsto \rho X_1$,  $e_i \mapsto E_i$,  and $g_i \mapsto G_i$ for
$1 \le i \le n-1$.
\end{theorem}

\subsection{The cyclotomic BMW and Kauffman tangle algebras}
 
We will consider algebras defined over ground rings with parameters satisfying certain conditions.  In order to avoid repeating these conditions, we establish the following convention:

\begin{definition}  A {\em ground ring} $S$ is a commutative unital ring with parameters $\rho$, $q$, $\delta_j$   ($j \ge 0$),  and
$u_1, \dots, u_r$, with    $\rho$, $q$,  $\delta_0$, and $u_1, \dots, u_r$  invertible, and with $\rho\inv - \rho=   (q\inv -q) (\delta_0 - 1)$.
\end{definition}

 \begin{definition} \label{definition:  cyclotomic BMW}
Let $S$ be a ground ring with
parameters $\rho$, $q$, $\delta_j$   ($j \ge 0$),  and
$u_1, \dots, u_r$.
The {\em cyclotomic BMW algebra}  $\bmw{n, S, r}(u_1, \dots, u_r)$
is the quotient of $\abmw{n, S}$ by the relation
\begin{equation} \label{equation: cyclotomic relation1}
(y_1 - u_1)(y_1 - u_2) \cdots (y_1 - u_r) = 0.
\end{equation}
\end{definition}

\begin{remark}
The assignment $e_i \mapsto e_i$, $g_i \mapsto g_i$, 
$y_1 \mapsto y_1$ defines a homomorphism $\iota$ from $\bmw  {n, S,r}$ to $\bmw  
{n+1,
S,r }$, since the relations are preserved.   It is not evident that $\iota$ is 
injective.  However, when $S$ is an admissible integral domain,  we will show that $\bmw  {n, S,r}$ is
 isomorphic to  a cyclotomic version of the Kauffman tangle algebra (defined  below), and in this case, it is true that $\iota$ is injective.
\end{remark}

\begin{remark} \label{remark: change of base ring for cyclotomic BMW algebras}
Let $S$ be a ring with parameters $\rho$, $q$, etc., as above, and let
 $S'$  be another  ring with parameters $\rho'$, $q'$, etc.  
 Suppose there is a  ring homomorphism $\psi : S \rightarrow S'$ mapping $\rho \mapsto \rho'$,  $q \mapsto q'$,  etc.  Then $\bmw{n, S', r}$  can be regarded as an algebra over $S$, with $s x = \psi(s) x$  for $s \in S$ and $x \in \bmw{n, S', r}$, and there is an
 $S$--algebra homomorphism $\tilde\psi : \bmw{n, S, r} \to \bmw{n, S', r}$  taking generators to generators.  We claim that $\bmw{n, S, r} \otimes_S S' \cong \bmw{n, S', r}$ as $S'$--algebras.  In fact, we have the $S'$ algebra homomorphism
$\tilde\psi \otimes \id: \bmw{n, S, r} \otimes_S S'  \to \bmw{n, S', r} \otimes_S S' \cong
\bmw{n, S', r}$.  In the other direction, we have an $S'$--algebra homomorphism
$\theta: \bmw{n, S', r} \to  \bmw{n, S, r} \otimes_S S'$  mapping $g_i \mapsto g_i \otimes 1$,  etc.  The maps $\tilde \psi \otimes \id$ and $\theta$  are inverses.

Of course,  this remark applies in general to algebras defined by generators and relations.\footnote{ We could have slightly  simplified some arguments in ~\cite{GH1} using this remark.}
\end{remark}

 Now we consider how to define a 
cyclotomic version of the Kauffman tangle algebra.
Rewrite the relation (\ref{equation: cyclotomic relation1})  in the form
\begin{equation*} \label{equation: cyclotomic relations2} 
\sum_{k = 0}^r  (-1)^{r-k}  \varepsilon_{r-k}(u_1, \dots, u_r) y_1^k = 0,
\end{equation*}
where $\varepsilon_j$  is the $j$--th elementary symmetric function.  The corresponding relation in the affine Kauffman tangle algebra is
\begin{equation} \label{equation: cyclotomic relations3} 
\sum_{k = 0}^r  (-1)^{r-k}  \varepsilon_{r-k}(u_1, \dots, u_r) \rho^k X_1^k = 0,
\end{equation}
Now we want to impose this as a local skein relation.

 \begin{definition}
Let $S$ be a ground ring with
parameters $\rho$, $q$, $\delta_j$   ($j \ge 0$),  and
$u_1, \dots, u_r$. 
The {\em cyclotomic Kauffman tangle algebra} $\kt{n, S, r}(u_1,  \dots, u_r)$ is the quotient of the affine Kauffman tangle algebra $\akt{n, S}$ by the
cyclotomic skein relation:
\begin{equation} \label{equation: kt cyclotomic relation}
\sum_{k = 0}^r  (-1)^{r-k}  \varepsilon_{r-k}(u_1, \dots, u_r) \rho^k \inlinegraphic[scale=.7]{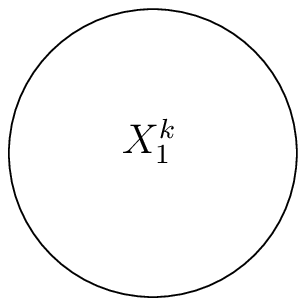} = 0,
\end{equation}
The sum is over affine tangle diagrams which differ only in the interior of  the indicated disc and
are identical outside of the disc;    the interior of the disc contains an interval on the flagpole and a piece of an affine tangle diagram isotopic to $X_1^k$.
\end{definition}

We continue to write $E_i$,  $G_i$, $X_1$   for the image of these elements
of the affine Kauffman tangle algebra in the cyclotomic Kauffman tangle algebra.
We write
$$Y_1 = \rho X_1.$$

The ideal in the affine Kauffman tangle algebra by which we are taking the quotient contains the ideal generated by
$
(Y_1 - u_1)(Y_1 - u_2) \cdots (Y_1 - u_r),
$
but could in principal be larger.  
It follows that there is a homomorphism  $$\varphi: \bmw{n, S, r}(u_1, \dots, u_r) \to \kt{n, S, r}(u_1, \dots, u_r)$$  determined by $y_1 \mapsto   Y_1$,   $e_i \mapsto E_i$,  $g_i  \mapsto G_i$.  Moreover, the following diagram (in which the vertical arrows are the quotient maps) commutes.    
\begin{diagram} \label{diagram: hm from cyclotomic bmw to cyclotomic kt}
\abmw {n, S }    ¤\Ear {\scriptstyle\varphi}     ¤\akt {n, S} ¤¤
\Sar {\scriptstyle\pi}         ¤          ¤\Sar {\scriptstyle\pi}  ¤¤
\movevertexleft{\bmw {n, S, r}(u_1, \dots, u_r)}                 ¤\Ear {\scriptstyle\varphi}     ¤\movevertexright{\kt {n, S, r}(u_1, \dots, u_r)}  ¤¤
\end{diagram}

\noindent
The homomorphism  $\varphi: \bmw{n, S, r}(u_1, \dots, u_r) \to \kt{n, S, r}(u_1, \dots, u_r)$  is surjective  because  the diagram commutes and
$\varphi: \abmw{n, S} \to \akt{n, S}$ is an isomorphism.

\begin{remark} \label{remark: change of base ring for tangle algebras}
 Let $S$ be a ring with parameters $\rho$, $q$, etc., as above, and let
 $S'$  be another  ring with parameters $\rho'$, $q'$, etc.  
 Suppose there is a  ring homomorphism $\psi : S \rightarrow S'$ mapping $\rho \mapsto \rho'$,  $q \mapsto q'$,  etc.    Any $S'$--algebra can be regarded as an $S$--algebra using $\psi$.
 We have a map of monoid rings
$
\tilde \psi : S \ \uhat n n  \rightarrow S' \ \uhat n n ,
$
and this map respects regular isotopy,  the Kauffman skein relations, and the cyclotomic relations, so induces an $S$--algebra homomorphism from $\tilde \psi:  \kt{n, S, r} \to \kt{n, S', r}$.   As in Remark \ref{remark: change of base ring for cyclotomic BMW algebras}, we have $\kt{n,S,r} \otimes_S S' \cong \kt{n,S', r}$ as $S'$--algebras.
\end{remark}

\subsection{Morphisms of ground rings}

 We consider what are the appropriate morphisms between ground rings for cyclotomic BMW or Kauffman tangle algebras.   The obvious notion would be that of a ring homomorphism taking parameters to parameters;  that is,  if $S$ is a ring with parameters $\rho$, $q$, 
etc.,  and $S'$ a weakly admissible ring with parameters $\rho'$, $q'$, etc.,  then a morphism $\varphi : S \to S'$  would be required to map $\rho \mapsto \rho'$,  
$q \mapsto q'$,  etc.  

However,  it is better to require less, for the following reason:  The parameter $q$ enters into the cyclotomic BMW
or Kauffman tangle relations only in the expression $q\inv -q$, and the transformation $q \mapsto -q\inv$ leaves this expression invariant.  Moreover,  the transformation $g_i  \mapsto -g_i$,  $\rho \mapsto -\rho$,  $q \mapsto -q$ (with all other generators and parameters unchanged)   leaves the cyclotomic BMW relations unchanged.   
Likewise, for the cyclotomic Kauffan tangle algebras,  the map taking an affine tangle diagram $T$ to $(-1)^{\kappa(T)} T$,  where $\kappa(T)$ is the number of crossings of ordinary strands  of $T$,   determines an isomorphism   $\kt {n, S, r}(\rho, q, \dots) \cong
\kt {n, S, r}(-\rho, -q, \dots)$.

Taking this into account, we arrive at the following notion:

\begin{definition}  \label{definition: parameter preserving}
Let $S$ be a ground ring with
parameters $\rho$, $q$, $\delta_j$   ($j \ge 0$),  and
$u_1, \dots, u_r$.
Let  $S'$ be another ground   ring with parameters $\rho'$, $q'$, etc.  

A unital ring homomorphism $\varphi : S \rightarrow S'$ is a {\em morphism of ground rings}  if it maps
$$
\begin{cases}
&\rho \mapsto \rho',  \text{ and}\\
& q \mapsto q' \text{ or } q \mapsto    -{q'}\inv,
\end{cases}
$$
or
$$
\begin{cases}
&\rho \mapsto  -\rho',  \text{ and}\\
& q \mapsto -q' \text{ or } q \mapsto    {q'}\inv,
\end{cases}
$$
and strictly preserves all other parameters.
\end{definition}

 Suppose there is a morphism of ground rings $\psi : S \rightarrow S'$.
  Then $\psi$ extends to a 
homomorphism  from 
$\bmw{n, S, r}$  to $\bmw{n, S', r}$,  and likewise to a homomorphism from
$\kt{n, S, r} $ to  $ \kt{n, S', r}$.     Moreover,  $\bmw{n, S, r} \otimes_S S' \cong \bmw{n,S', r}$  and $\kt{n, S, r} \otimes_S S' \cong \kt{n, S', r}$ as $S'$--algebras.  This follows from the observations above and Remarks \ref{remark: change of base ring for cyclotomic BMW algebras}
  and \ref{remark: change of base ring for tangle algebras}.

\subsection{Weak admissibility}  \label{subsection: weak admissibility}
As noted in 
~\cite{GH2}, Section 4.1,
unless the parameters of the ground ring $S$
satisfy the {\em weak admissibility} condition defined below,   the identity element $\bm 1$ of the cyclotomic Kauffman tangle algebras
will be a torsion element over $S$;  if $S$ is a field (and the additional relations do not hold)  then $\bm 1 = 0$,  so $\kt {n, S, r} = \{0\}$.
Similarly,  in the absence of weak admissibiltity,  a computation done in the cyclotomic BMW algebra shows that  $e_1 \in \bmw {2, S, r}$ 
  is a torsion element over $S$.

\begin{definition} \label{definition: weak admissibility}
Let $S$ be a ground ring with parameters
$\rho$, $q$,  $\delta_j$, $j \ge 0$, and $u_1, \dots, u_r$.
  We say that the parameters are {\em weakly admissible}  (or that the ring $S$ is weakly
admissible)  if the  following relation holds:
$$
\sum_{k = 0}^r  (-1)^{r-k}  \varepsilon_{r-k}(u_1, \dots, u_r)  \delta_{k + a}  = 0,
$$
for $a \in \Z$,  where for $j \ge 1$,   $\delta_{-j}$ is defined by  the recursive relations of Lemma \ref{lemma - recursion for f_r}.
\end{definition}

\subsection{Inclusion, conditional expectation, and trace} \label{subsection: Inclusion, conditional expectation, and trace}
Let $S$  be a ground ring with  parameters as above and write
$\kt{n, S, r}$  for $\kt{n, S, r}(u_1, \dots, u_r)$.
The affine Kauffman tangle algebras have inclusion maps
 $\iota : \akt {n-1} \rightarrow \akt {n}$, defined on the level of affine tangle diagrams
by adding an additional strand on the right
without adding any crossings:
$$
\iota: \quad \inlinegraphic{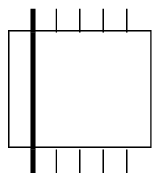} \quad \mapsto \quad 
\inlinegraphic{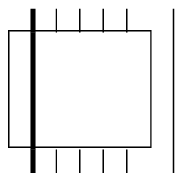}.
$$
Since these maps respect the cyclotomic relation (\ref{equation: kt cyclotomic relation}),  they induce homomorphisms
$$\iota : \kt {n-1, S, r} \rightarrow \kt {n, S, r}.$$
Recall also  that the affine Kauffman tangle algebras have a conditional expectation 
$\eps_n : \akt n \rightarrow \akt {n-1}$ defined by
$$
\eps_n(T) = \delta_0\inv{\rm cl}_n(T), 
$$
where   ${\rm cl}_n$ is the map of affine $(n,n)$--tangle diagrams to affine
$(n-1, n-1)$--tangle diagrams that ``closes" the rightmost strand, 
without adding any crossings:
$$
{\rm cl}_n: \quad \inlinegraphic{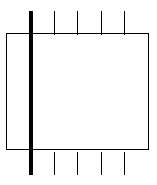} \quad \mapsto \quad 
\inlinegraphic{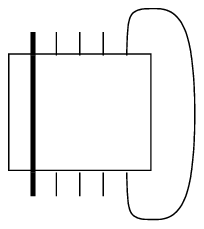}.
$$
These maps respect the cyclotomic relation (\ref{equation: kt cyclotomic relation}),  so 
induce conditional expectations
$$\eps_n : \kt{n, S, r} \rightarrow \kt{n-1, S, r}.$$  
Since $\eps_n\circ \iota$ is the identity on $\kt{n-1, S, r}$,  it follows that 
$\iota : \kt {n-1, S, r} \rightarrow \kt {n, S, r}$
is injective.

It is shown in the companion paper ~\cite{GH2},  Proposition 4.3
 that if the ground ring
$S$ is weakly admissible, then the cyclotomic Kauffman tangle algebra $\kt{0, S, r}$ is a free $S$--module of rank 1.  Define 
$$
\eps = \eps_{1}\circ\cdots \circ \eps_n  : \kt{n, S, r} \to \kt{0, S, r} \cong S
$$
It follows from ~\cite{GH1}, Proposition 2.14, that $\eps$ is a trace.
We also define $\eps: \bmw{n, S, r} \to S$ by $\eps = \eps \circ \varphi$, 
where $\varphi: \bmw{n, S, r} \to \kt{n, S, r}$ is the canonical homomorphism.  Then
$\eps$ is a trace on $\bmw{n, S, R}$ with the {\em Markov property}:
 for $b \in \bmw {n-1, S, r}$,
\begin{enumerate}
\item[\rm(a)] $\eps(b g_{n-1}^{\pm 1} ) = (\rho^{\pm 1}/\delta_0) \eps(b) $,
\item[\rm(b)] $\eps(b e_{n-1}) = (1/\delta_0) \eps(b)$,
\item[\rm(c)\hskip1.2pt]  $\eps(b (y'_n)^r) = \delta_{r} \eps(b)$, and
\end{enumerate}
for $r \in \Z$, 
where $y'_n = (g_{n-1} \cdots g_1) y_1 (g_1\inv \cdots g_{n-1}\inv)$.  See ~\cite{GH1}, Corollary 6.16.

\begin{lemma} \label{lemma:  formulas involving cond exp and En}  Let $S$ be a weakly admissible ring. 
For $x \in \kt {n, S, r}$,
\begin{enumerate}
\item  $E_n x E_n = \delta_0  \eps_n(x)$.
\item  $\eps(E_n x E_n) = \eps(x)$.
\item  $\eps(E_n x) =  \delta_0\inv  \eps(x)$.
\end{enumerate}
\end{lemma}

\begin{proof}  Straightforward.
\end{proof}

\begin{lemma} \label{lemma:  E and G not zero}
Let $S$ be a weakly admissible ring.  Then for all $n \ge 1$, $E_n$ and  $G_n$ are non-zero elements in $\kt {n, S, r}$.
\end{lemma}

\begin{proof}
$\eps(E_n) = \delta_0\inv$  and $\eps(G_n) = \rho\inv \delta_0^{-2}$.
\end{proof}

\subsection{Inclusions of split semisimple algebras}
 A general source for the material in this section is ~\cite{GHJ}.   The Jones basic construction (discussed below) is from ~\cite{Jones}.
 
 A finite dimensional split semisimple algebra over a field $F$ is one which is isomorphic to a finite direct sum of full matrix algebras over $F$.
Suppose
$A$ and $B$ are finite dimensional split semisimple algebras over  $F$.  Let   $A_i$,  $i \in I$, be the minimal ideals of $A$  and  $B_j$,  $j \in J$,   the minimal ideals of $B$.  Let 
$\psi : A \rightarrow B$ be a unital homomorphism of algebras.   

We associated a  $J \times I$   {\em multiplicity matrix}
$\Lambda_\psi$ to $\psi$, as follows.  Let $W_j$ be a simple $B_j$--module.
Then $W_j$ becomes an $A$--module via $\psi$,  and $\Lambda_\psi(j, i)$ is the multiplicity of  a simple $A_i$--module 
 in the decomposition of $W_j$ as an $A$--module.   In particular,  if $A \subseteq B$ (and both algebras have the same identity element)  the multiplicity matrix for the inclusion is called  the {\em inclusion matrix} for $A \subseteq B$.   

An equivalent characterization of the multiplicity matrix $\Lambda_\psi$ of $\psi$ is the following.   Let $q_i$ be a minimal idempotent in $A_i$  and let $z_j$  be the identity of 
$B_j$  (a minimal central idempotent in $B$).  Then $\psi(q_i) z_j$ is the sum of
 $\Lambda_\psi(j, i)$ minimal idempotents in $B_j$.

Sometimes it is convenient to encode an inclusion matrix (or multiplicity matrix)  by a bipartite graph, called the {\em branching diagram};  the branching diagram has vertices labeled by $I$  arranged on one horizontal line,  vertices labelled by  $J$  arranged along a second (higher) horizontal line,    and $\Lambda_\psi(j, i)$ edges connecting
$j \in J$ to $i \in I$.

If $A \spp 1 \subseteq A \spp 2 \subseteq A \spp 3  \cdots$ is a sequence of inclusions of finite dimensional split semisimple algebras over $F$,  then the branching diagram for the sequence is obtained by stacking the branching diagrams for each inclusion,  with the 
upper vertices of the diagram for $A \spp i \subseteq A \spp {i+1}$  being identified with the lower vertices of the diagram for $A \spp {i+1} \subseteq A \spp {i+2}$.

Let $A \subseteq B$ be a pair of finite dimensional split semisimple algebras over $F$.
The {\em Jones basic construction}  for $A \subseteq B$ is the pair
$B \subseteq  \End(B_A)$.    That is,  $B$ is regarded as acting on itself by left multiplication,  and is thus included in the endomorphisms of $B$ as a right $A$--module.
Then $\End(B_A)$ is again a split semisimple algebra whose minimal ideals are in one--to--one correspondence with those of $A$,  and the inclusion matrix for the pair $B \subseteq  \End(B_A)$ is the transpose of that for $A \subseteq B$.

Suppose now that $B$ has a faithful trace $\eps$ with faithful restriction to $A$.
Then there is a unique trace preserving conditional expectation (i.e.  $A$--$A$ bimodule map)   $\eps_A : B \to A$  determined by $\eps (b a) =  \eps( \eps_A(b) a)$  for all
$a \in A$ and $b\in B$.  Note that $\eps_A$ is an idempotent in $\End(B_A)$.  The Jones basic construction has a particular realization in terms of this data. In fact,
$\End(B_A)  =  \langle B, \eps_A \rangle =  B \eps_A  B,$  
where $ \langle B, \eps_A \rangle $ denotes the subalagebra of $\End(B_A)$ generated
by $B$ and $\eps_A$.
Moreover,  if $q_i$ is a minimal idempotent in $A_i$,  then
$q_i  \eps_A = \eps_A q_i$ is a minimal idempotent in the corresponding minimal ideal of 
$\End(B_A)$.

\begin{theorem}[Wenzl \cite{Wenzl-Brauer}] \label{Theorem: Wenzl extensions} Suppose $A \subseteq B \subseteq C$ are $F$--algebras with
$A$ and $B$ finite dimensional split semisimple.   Suppose that $B$ has a trace $\eps \in B^*$ that is non-degenerate on $B$ and has non-degenerate restriction to $A$.  
Let $\eps_A$  denote the corresponding 
trace preserving conditional expectation
 $\eps_A : B \rightarrow A$.

Suppose that $e \in C$ is an idempotent such that $e x e = \eps_A(x)  e$  for all $x \in B$ and
$a \mapsto a e$ is an injective homomorphism from $A$ to $C$.   Let $\langle B, e \rangle$ be the subalgebra of $C$ generated by $B$ and $e$.    Then
$B e B$ is an ideal in $\langle B, e \rangle$, and $$B e B \cong B \eps_A B = \langle B, \eps_A\rangle = \End(B_A).$$   
\end{theorem}

\begin{remark} \label{remark:  on Wenzl extensions}
 In the situation of the theorem,  it follows that $BeB$  is split semisimple, with minimal ideals in one--to--one correspondence with those of $A$.   If $q_i$ is a minimal idempotent in a minimal ideal $A_i$ of $A$,  then $q_i e$ is a minimal idempotent in the corresponding minimal ideal of $B e B$.    $BeB$ has an identity element $z$,  necessarily a central idempotent in $\langle B, e \rangle$.   The multiplicity matrix for the homomorphism $b \mapsto z b$  from $B$ to $B e B$ is the transpose of the inclusion matrix for $A \subseteq B$.
\end{remark}

We now discuss {\em path idempotents} for a sequence $A \spp 0 \subseteq A \spp 1 \subseteq A \spp 2 \cdots$ of split semisimple algebras, under the simplifying assumption  (sufficient for our purposes)  that $A \spp 0 = F$ and each inclusion $A \spp i \subseteq A\spp {i+1}$ is {\em multiplicity--free}, that is, the branching diagram has at most one edge connecting any pair of vertices.  Let $\Gamma_k$  be an index set for the minimal ideals of 
$A\spp k$;  denote the minimal ideals of $A\spp k$ by $A\spp k_\lambda$ and the minimal central idempotents by  $z\spp k_\lambda$  ($\lambda \in \Gamma_k$).   
Denote the unique element of $\Gamma_0$ by $\emptyset$.
For 
$\mu \in \Gamma_j$ and $\lambda \in \Gamma_{j+1}$,  write
$\mu \subseteq \la$ if $\mu$ and $\lambda$ are connected by an edge of the branching diagram.  A {\em path} of length $k$ in the branching diagram is a sequence
$(\emptyset, \lambda\spp 1, \lambda \spp 2, \dots, \lambda \spp k)$  with
$\lambda \spp j \subseteq \lambda \spp {j+1}$  for all $j$.  Let $\mathcal T(k)$  be the set of all paths of length $k$,  and $\mathcal T(k, \lambda)$  the set of paths of length $k$ with final vertex $\lambda$.   If $T = (\emptyset, \lambda \spp 1, \dots, \lambda \spp k)$
is a path of length $k$,  let $T'$ denote its truncation of length $k-1$,
$T' = (\emptyset, \lambda \spp 1, \dots, \lambda \spp {k-1})$.  Moreover, if $\lambda \in \Gamma_{k+1}$ and
$\lambda \spp k \subseteq \lambda$,  let $T + \lambda$  denote the extension of $T$ of length
$k+1$ with final vertex $\lambda$.
We claim that for each $k$, there is a canonical family of idempotents
$\{p_T : T \in \mathcal T(k)\}$, in $A \spp k$, with the properties:
\begin{enumerate}
\item  The elements $p_T$ are mutually orthogonal minimal idempotents in $A \spp k$, and the sum of the $p_T$ with
$T \in \mathcal T(k, \lambda)$ is the minimal central idempotent $z \spp k_\lambda$.
\item  $p_T p_{T'} = p_{T'} p_T = p_T$.
\end{enumerate}
The construction is simple:  Let $\emptyset$ denote the path of length $0$.  Necessarily
$p_\emptyset = 1$.   If the idempotents $p_T$  have been constructed for paths of length $j$,  $T \in \mathcal T(j, \mu)$,  and $\lambda \in \Gamma_{j+1}$ satisfies
$\mu \subseteq \lambda$, then put $p_{T + \lambda}  =  p_T z \spp {j+1}_\lambda$, which is a minimal idempotent in $A \spp {j+1}$  because of the assumption of multiplicity--free inclusions.   

It follows that if $T = (\emptyset, \lambda\spp 1, \lambda \spp 2, \dots, \lambda \spp k)$ is a path of length $k$,  then for all $j \le k$,  $A \spp j p_T =  A\spp j_{\lambda \spp j} p_T$  affords the simple $ A\spp j_{\lambda \spp j}$--module.

\begin{lemma} \label{lemma:  conditional expectation of path idempotents}
 Suppose $A \spp 0 \subseteq A \spp 1 \subseteq A \spp 2 \subseteq \cdots \subseteq A \spp n$  is a sequence of split semisimple algebras, with  $A \spp 0 = F$ and with each inclusion $A \spp i \subseteq A\spp {i+1}$  { multiplicity--free}.
Suppose we have a trace $\eps$ on $A_n$,  not necessarily faithful,  but with faithful restriction to $A_{n-1}$.  Suppose also that we have a trace--preserving conditional expectation $\eps_n: A_n \to A_{n-1}$.  
Let $T$ be a path of length $n$ on the branching diagram for the sequence of inclusions.   Then 
\begin{equation*}
 \eps_{n}(p_T) =   (\eps(p_T)/ \eps(p_{T'})) p_{T'}.
 \end{equation*}
\end{lemma}

\begin{proof}
We have  
 $
 \eps_{n}(p_T)  =   \eps_{n}(p_{T'} p_T p_{T'}) = p_{T'} \eps_{n}(p_T) p_{T'}.
 $
 Since $p_{T'}$ is a minimal idempotent in $A_{n-1}$,  we have $ \eps_{n}(p_T) = \kappa p_{T'}$
 for some $\kappa \in F$.    Taking traces,  we have  $\eps(p_T) = \kappa  \eps(p_{T'})$, so
 $\kappa =  \eps(p_T)/ \eps(p_{T'})$.
\end{proof}

\section{Representations of the 2--strand algebra and admissibility conditions}
\label{section:  representations and admissibility}
To obtain any substantial results about the cyclotomic BMW  and Kauffman tangle algebras, it appears to be necessary to impose a stronger  condition on the ground ring than weak admissibility.  Appropriate conditions can be found by considering the representation theory of the the 2--strand algebra, in particular  representations in which the generator $e_1$ acts non--trivially.   

This idea was used by Ariki, Mathas, and Rui in a related context ~\cite{ariki-mathas-rui};  the analogue of their condition,  translated to our context, is useful, but too restrictive.
The ``correct" notion of admissibility for cyclotomic BMW algebras was discovered
 by Wilcox and Yu  
 ~\cite{Wilcox-Yu}.
The left module $\bmw {2, S, r} \, e_1$  is spanned over the ground ring by $\{e_1, y_1 e_1, \dots,  y_1^{r-1} e_1\}$, as is easily checked.  Generically, one would expect this spanning set to be linearly independent over $S$, and the requirement of linear independence is the module--theoretic formulation of the admissibility condition of Wilcox and Yu.

\subsection{Admissibility} \label{subsection: admissibility}

Consider a ground  ring $S$ with parameters $\rho$, $q$, $\delta_j$ ($j \ge 0$)  and $u_1, \dots, u_r$. 
Let $a_j$ denote the signed elementary symmetric function in $u_1, \dots, u_r$, 
$
a_j = (-1)^{r-j} \varepsilon_{r-j}(u_1, \dots, u_r).  
$
Let
 $W_2$  denote the cyclotomic BMW algebra
$W_2 = \bmw {2, S, r}(u_1, \dots, u_r)$.  Write $e$ for $e_1$ and $g$ for $g_1$.

\begin{lemma} \label{lemma: powers of y generate W2e}   The left ideal $W_2\, e$ in $W_2$  is equal to the $S$--span of  $\{e, y_1 e, \dots, y_1^{r-1} e\}$.
\end{lemma}

\begin{theorem}[Wilcox-Yu, \cite{Wilcox-Yu}] \label{theorem: equivalent conditions for admissibility}
 Let $S$ be a ground ring with
parameters $\rho$, $q$, $\delta_j$ ($j \ge 0$)  and $u_1, \dots, u_r$.  
Assume that $(q - q\inv)$ is non--zero and not a zero--divisor in $S$.
The following conditions are equivalent:
\begin{enumerate}
\item  $S$ is weakly admissible, and $\{e, y_1 e, \dots, y_1^{r-1} e\} \subseteq W_2$ is linearly independent over $S$.
\item  The parameters satisfy the following relations:

\begin{equation} \label{equation: yu wilcox admissibility condition 1}
\begin{aligned}
&\rho(a_\ell - a_{r-\ell}/a_0) \ + 
\\& (q-q\inv)\bigg [ \sum_{j = 1}^{r - \ell} a_{j+\ell} \delta_j 
-  \sum_{j = \max(\ell + 1, \lceil r/2 \rceil)}^{\lfloor (\ell + r)/2 \rfloor} a_{2j - \ell}
+   \sum_{j =  \lceil \ell/2 \rceil}^{\min(\ell, \lceil r/2 \rceil -1)} a_{2j - \ell} \bigg ]= 0, \\& \quad \text{for $1 \le \ell \le r-1$},
\end{aligned}
\end{equation}
\begin{equation} \label{equation: yu wilcox admissibility condition 2}
\rho\inv a_0 - \rho a_0\inv = 
\begin{cases}
0 & \text{if  $r$ is odd} \\
(q - q\inv) & \text{if $r$ is even},
\end{cases}
\end{equation}
 and
\begin{equation} \label{equation:  weak admissibility in yu wilcox}
\delta_a = -\sum_{j=0}^{r-1} a_j \delta_{a-r+j} \quad \text{for }  a \ge r.
\end{equation}

\item  $S$ is weakly admissible, and $W_2$ admits a module $M$ with an $S$--basis $\{v_0, y_1 v_0, \dots,  y_1^{r-1} v_0\}$  such that
$e v_0 = \delta_0 v_0$.    
\end{enumerate}
\end{theorem}

\begin{definition}[Wilcox and Yu, \cite{Wilcox-Yu}]     Let $S$ be a  ground ring with
parameters $\rho$, $q$, $\delta_j$ ($j \ge 0$)  and $u_1, \dots, u_r$.   One says that $S$ is {\em admissible} (or that the parameters are {\em admissible})   if $(q - q\inv)$ is non--zero and not a zero--divisor in $S$ and if the equivalent conditions of Theorem \ref{theorem: equivalent conditions for admissibility} hold.
\end{definition}

In the following statement, 
we  omit the assumption that
$q - q\inv$ is not a zero--divisor,  but we  impose the assumption that $S$ is weakly admissible. The statement can be proved by a minor modification of the proof of Wilcox and Yu.

\begin{corollary}  \label{corollary:  equivalent conditions for "admissiblity" without zero divisor assumption}
Let $S$ be a weakly admissible ring with parameters
$\rho$, $q$,  $\delta_j$  ($j \ge 0$), and $u_1, \dots, u_r$.  The following conditions are equivalent.
\begin{enumerate}
\item  $\{e, y_1 e, \dots, y_1^{r-1} e\}  \subseteq W_2$ is linearly independent over $S$.
\item  The parameters satisfy the relations (\ref{equation: yu wilcox admissibility condition 1}) and (\ref{equation: yu wilcox admissibility condition 2}) of Wilcox and Yu.

\item  $W_2$ admits a cyclic module $M$ with an $S$--basis  $\{v_0, y_1 v_0, \dots,  y_1^{r-1} v_0\}$  such that
$e v_0 = \delta_0 v_0$.    
\end{enumerate}
\end{corollary}

\vbox{
\begin{remark} \label{remark: parameter preserving transformations preserve admissibility}
 \mbox{}
\begin{enumerate}
\item
Let $S$ be a ground ring with admissible parameters  $\rho$, $q$, $\delta_j$   ($j \ge 0$),  and
$u_1, \dots, u_r$.  Then  
$$
\rho, -q\inv, \delta_j \  (j \ge 0),  \text{ and }  u_1, \dots, u_r
$$
and
$$
-\rho, -q, \delta_j \  (j \ge 0),  \text{ and }  u_1, \dots, u_r
$$
are also sets of admissible parameters.  This follows from condition (2) of  Theorem \ref{theorem: equivalent conditions for admissibility}.  
\item  If $S$ is a  ground ring with admissible parameters and $\varphi : S \rightarrow S'$ is 
a morphism of ground rings in the sense of Definition \ref{definition: parameter preserving}, such that
\break $\varphi(q - q\inv)$  is non--zero and not a zero divisor,  then $S'$ is also admissible.
\end{enumerate}
\end{remark}
}

\begin{remark}  \label{remark:  solving for rho and deltas}
Assume that $S$ is a ground  ring with admissible parameters.
We observe that the admissibility relation (\ref{equation: yu wilcox admissibility condition 1}) can be solved for 
 $$(q - q\inv) \delta_1,  \dots,  (q - q\inv) \delta_{r-1}$$ in terms of  $\rho$, 
$a_0\powerpm, a_1, \dots, a_{r-1}$, and $q^{\pm 1}$.  
Rewrite the relations
(\ref{equation: yu wilcox admissibility condition 1}) in the form
\begin{equation}
(q - q\inv)  \sum_{j=1}^{r-\ell} a_{j+\ell}  \delta_j  = p_\ell(q^{\pm 1}, \rho, a_0\powerpm, a_1, \dots, a_{r-1}) \quad \text(1 \le \ell \le r-1),
\end{equation}
where the $p_\ell$ are specific polynomials.   This is a triangular system of linear equations in the variables $(q - q\inv) \delta_j$  with a unitriangular matrix of coefficients
$$
\begin{bmatrix}
1 & & & &\\
a_{r-1} & 1 & && \\
a_{r-2} & a_{r-1} & 1 & &\\
\vdots  & & \ddots& \ddots & \\
a_2   &a_3  &\dots&a_{r-1}& 1\\
\end{bmatrix},
$$
so it has a solution of the form
$$
(q - q\inv) \delta_j =  Q_j(q^{\pm 1}, \rho,  a_0\powerpm , a_1, \dots, a_{r-1}), \quad \text{for} \ \  1 \le j \le r-1,
$$
where the $Q_j$ are again some specific polynomials.  Because $q - q\inv$ is assumed not to be a zero--divisor in $S$,  we can append an inverse of $q - q\inv$ to $S$.  In
$S[(q - q\inv)\inv]$, we have
$$
 \delta_j = (q - q\inv)\inv Q_j(q^{\pm 1}, \rho,  a_0\powerpm , a_1, \dots, a_{r-1}), \quad \text{for} \ \ 1 \le j \le r-1,
$$
\end{remark}

\begin{remark} \label{remark:  solve for rho if S is an integral domain}
Equation (\ref{equation: yu wilcox admissibility condition 2}) is equivalent to
$$
\begin{cases}
(a_0 - \rho)(a_0 + \rho) = 0 &\text{if $r$ is odd, and}\\
(\rho - q\inv a_0)(\rho + q a_0) = 0 &\text{if $r$ is even.}\\
\end{cases}
$$
If $S$ is an integral domain, we have  $\rho = \pm a_0$ if $r$ is odd,  and $\rho = -q a_0$ or $\rho = q\inv a_0$ if $r$ is even.  However,  we do not need to entertain both solutions.  Considering Remark \ref{remark: parameter preserving transformations preserve admissibility},  we can assume
$$\rho = -a_0 =\prod_{j = 1}^r u_j,$$
if $r$ is odd, and
$$\rho = q\inv a_0 = q\inv \prod_{j = 1}^r u_j,$$
if $r$ is even.
\end{remark}

We have the following consequence of this discussion:

\ignore{
\begin{corollary} \label{corollary:  rho and deltas contained in subring}
Let $S$ be a ring with admissible parameters $\rho$, $q$, etc.  
\begin{enumerate}
\item
  The elements $\delta_j$  for $j \ge 0$  are contained in  the subring of
$S[(q-q\inv)\inv]$ 
generated by $\rho, q^{\pm 1}$, $(q\inv - q)\inv$, and 
$u\powerpm_1, \dots, u\powerpm_r$.
\item  If we can write $\rho$ as a monomial in $q\powerpm$ and $u_1, \dots, u_r$, then
  then $\rho\powerpm$  and $\delta_j$  for $j \ge 0$  are elements of the subring of
$S[(q-q\inv)\inv]$ 
generated by $q^{\pm 1}$, $(q\inv - q)\inv$, and 
$u\powerpm_1, \dots, u\powerpm_r$. This is so, in particular, if $S$ is an integral domain.
\end{enumerate}
\end{corollary}
}

\begin{corollary} \label{corollary:  rho and deltas contained in subring}
Let $S$ be an integral ground ring   with admissible parameters $\rho$, $q$, etc.  
Then $\rho\powerpm$  and $\delta_j$  for $j \ge 0$  are elements of the subring of
$S[(q-q\inv)\inv]$ 
generated by $q^{\pm 1}$, $(q\inv - q)\inv$, and 
$u\powerpm_1, \dots, u\powerpm_r$.
\end{corollary}

\subsection{$u$--admissibility} \label{subsection: u--admissibility}

 We maintain the notation established at the beginning of Section \ref{subsection: admissibility}.  In this section,  we assume that the ring $S$ is an integral domain,  that
 $q - q\inv \ne 0$,  that the parameters $u_i$  are distinct, and that $u_i u_j \ne 1$  for all $i, j$.

\begin{lemma} \label{corollary:  gammas defined by linear equations}
 Let $F$ be field and let $u_1, \dots, u_r$  be distinct non--zero elements of $F$ with $1 \ne u_i u_j$  for all $i, j$.   Let $\rho$ and $q$ be non--zero elements of $F$ with $q^2 \ne 1$.
Then the unique solution to the system of linear equations
\begin{equation}\label{equation:  equation for gammas00}
\sum_j  \frac{1}{1-u_i u_j}\,  \gamma_j =    \frac{1}{1-u_i^2}  +  \frac{1}{\rho (q\inv - q)} \quad (1 \le i \le r)
\end{equation}
is
\begin{equation}\label{equation: formula for gammas00}
\gamma_j =\prod_{\ell \ne j} \frac{(u_\ell u_j -1)}{u_j - u_\ell} \left(  
 \frac{1-u_j^2}{\rho(q\inv -q)}  \prod_{\ell \ne j} u_\ell   + 
\begin{cases}
1 &\text{ if $r$  is odd} \\
-u_j     &\text{if $r$ is even} 
\end{cases}
\right )
\end{equation}
\end{lemma}

We defer the proof of the lemma to  Remark 
\ref{proof of formula for gammas} in Section \ref{subsection:  symmetric functions}.

\begin{theorem}\label{theorem: equivalent conditions for u--admissibility}
 Let $S$ be an integral ground ring with parameters $\rho$, $q$, $\delta_j$ ($j \ge 0$)  and $u_1, \dots, u_r$.
    Assume that $(q - q\inv) \ne 0$, that   the elements $u_i$ are distinct,  and  that $u_i u_j \ne 1$ for all $i, j$.
 Define $\gamma_j$ in the field of fractions of $S$ by  (\ref{equation: formula for gammas00}), for $1 \le j \le r$.

The following conditions are equivalent:
\begin{enumerate}
\item  $S$ is weakly admissible, and $\{e, y_1 e, \dots, y_1^{r-1} e\} \subseteq  \bmw {2,S}$ is linearly independent over $S$.
\item    For all $a \ge 0$, we have $\delta_a = \sum_{j = 1}^r  \gamma_j u_j^a$.

\item  $S$ is weakly admissible, and $\bmw {2, S}$ admits a module $M$ with an $S$--basis $\{v_0, y_1 v_0, \dots,  y_1^{r-1} v_0\}$  such that
$e v_0 = \delta_0 v_0$.    
\end{enumerate}
\end{theorem}

\begin{proof}  Let $F$ denote the field of fractions of $S$,  and let (1$_F$)  and 
(3$_F$) indicate the analogue of conditions (1) and (3) over $F$, i.e.

(1$_F$) \ \ $F$ is weakly admissible, and $\{e, y_1 e, \dots, y_1^{r-1} e\} \subseteq  \bmw {2,F}$ is linearly independent over $F$.

(3$_F$)\ \  $F$ is weakly admissible, and $\bmw {2, F}$ admits a module $M$ with an $F$--basis $\{v_0, y_1 v_0, \dots,  y_1^{r-1} v_0\}$  such that
$e v_0 = \delta_0 v_0$.

According to Theorem \ref{theorem: equivalent conditions for admissibility},  conditions (1) and (3) are equivalent, and moreover,  they are equivalent to the second  condition  of Theorem \ref{theorem: equivalent conditions for admissibility}.  Since that condition is unchanged if we replace $S$ by $F$,  it follows that (1) and (3)  are equivalent to 
(1$_F$)  and  (3$_F$).  Therefore,  to prove the theorem, it suffices to show that
$$\text{(1$_F$)  $\implies$ (2) $\implies$   (3$_F$).}$$

Assume condition (1$_F$).   Recall that the span of $\{e, y_1 e, \dots, y_1^{r-1} e\} $ is equal to the left ideal $\bmw {2, F} \, e$.  Moreover,  $e y_1^j e = \delta_j e$,  so
$e \bmw {2, F} e = F e$.

For each $j$  ($1 \le j \le r$)   define $p_j$ by 
$$
p_j = \prod_{i \ne j} \frac{y_1 - u_i}{u_j - u_i} \in \bmw {1, F}.
$$
Then we have $y_1 p_j =  u_j p_j$,  $p_j^2 = p_j$,  and $\sum_j p_j = 1$.

Let $m_j = p_j e$.  Then  $m_j \ne 0$,  by the linear independence of  $\{e, y_1 e, \dots, y_1^{r-1} e\} $ over $F$,  $y_1 m_j = u_j m_j$,  and $\sum_j m_j = e$.    
The set $\{m_1, \dots, m_r\}$ is linearly independent, since the $m_j$ are eigenvectors for $y_1$ with distinct eigenvalues.  Because the dimension of $\bmw {2, F}\, e$ is $r$,
$\{m_1, \dots, m_r\}$ is a basis of $\bmw {2, F}\, e$.

Now we  define  $\kappa_j \in F$  by   $e m_j = \kappa_j e  = \kappa_j \sum_i m_i$   and $c_{i, j}$ by $g m_j = \sum_i c_{i,j} m_i$.
Note that $y_1 g (y_1^s e) = g\inv y_1\inv (y_1^s e)$,  by Remark \ref{remark: involution on affine BMW}.
Therefore, 
$$
y_1 g m_j = g\inv y_1\inv m_j = \left ( g + (q - q\inv)(e -1) \right ) y_1\inv m_j.
$$
We have
$$
y_1 g m_j = y_1 \sum_i c_{i, j } m_i = \sum_i u_i c_{i,j} m_i.
$$
On the other hand,
$$
\begin{aligned}
( g + & (q - q\inv)(e -1)  ) y_1\inv m_j  \\
&= u_j\inv 
\sum_i c_{i,j} m_i  +  u_j\inv  (q - q\inv) ( \kappa_j \sum_i m_i  -  m_j)
\end{aligned}
$$
Equating the coefficient of $m_i$ in the two expressions gives
$$
(u_i - u_j\inv) c_{i,j}=  (q - q\inv) u_j\inv (\kappa_j -  \delta(i, j)),
$$
where $\delta(i,j)$ denotes the Kronecker delta.
Solving, we get
\begin{equation}\label{equation: formula for cij}
c_{i,j} =  (q\inv - q)  \frac{\kappa_j - \delta(i,j)}{1 - u_i u_j}.
\end{equation}
If we apply both sides of the equation $g e = \rho\inv e$ to $m_k$,  we get
$$
\kappa_k \sum_{i, j} c_{i, j} m_i = \rho\inv \kappa_k \sum_i m_i.
$$
Choosing $k$ so that $\kappa_k \ne 0$, and equating the coefficients of $m_i$ on both sides gives
$
\sum_j c_{i, j} = \rho\inv
$
for all $i$.
Taking  (\ref{equation: formula for cij}) into account, we have for all $i$, 
\begin{equation}\label{equation:  equation for gammas}
\sum_j  \frac{1}{1-u_i u_j}\,  \kappa_j =    \frac{1}{1-u_i^2}  +  \frac{1}{\rho (q\inv - q)}.
\end{equation}
It follows from  (\ref{equation:  equation for gammas}) and 
Corollary  \ref{corollary:  gammas defined by linear equations}  that 
$\kappa_j$ is given by
\begin{equation}\label{equation: formula for gammas}
\kappa_j = \gamma_j =\prod_{\ell \ne j} \frac{(u_\ell u_j -1)}{u_j - u_\ell} \left(  
 \frac{1-u_j^2}{\rho(q\inv -q)}  \prod_{\ell \ne j} u_\ell   + 
\begin{cases}
1 &\text{ if $r$  is odd} \\
-u_j     &\text{if $r$ is even} 
\end{cases}
\right )
\end{equation}
Next,  note that for $a \ge 0$
$$
\begin{aligned}
\delta_a e  &= e y_1^a e  = e y_1^a \sum_i m_i =   e \sum_i  u_i^a  m_i =  (\sum_i \gamma_i u_i^a)  e.
\end{aligned}
$$
Therefore, for $a \ge 0$,  $ \delta_a =  \sum_i  \gamma_i  u_i^a$.
This completes the proof of the implication (1$_F$) $\implies$  (2).

We now assume condition (2) and prove (3$_F$).  Let $\{\unitvector e_1, \dots, \unitvector e_r\}$ denote the standard basis of 
$M = F^r$.
Define linear maps $E, G, Y$ on $M$ by
\begin{equation} \label{equation: formulas for r diml repn2}
\begin{split}
&Y \,\unitvector e_j = u_j \unitvector e_j, \\ 
&E\, \unitvector e_j =  \gamma_j (\unitvector e_1 +  \cdots + \unitvector e_r)\\ 
&G\, \unitvector e_j =   \displaystyle (q\inv - q)  \sum_{i = 1}^r  \frac{\gamma_j - \delta(i,j)}{1 - u_i u_j} \ \unitvector e_i,\\
\end{split}
\end{equation}
where for $1 \le j \le r$,  $\gamma_j$ is given by  (\ref{equation: formula for gammas00}).

Applying the definitions of $Y, G, E$, and the $u$--admissibility assumption, we have, for $a \ge 0$, 
$
E Y^a E\,  \unitvector e_j = ( \sum_i \gamma_i u_i^a) E \, \unitvector e_j = \delta_a E \, \unitvector e_j.
$
Thus,  for $a \ge 0$, 
\begin{equation} \label{equation:  BMW relations for E G Y 0}
E Y^a E = \delta_a E,
\end{equation}
and, in particular, $E^2 = \delta_0 E$.
It is straightforward to compute  that
\begin{equation} \label{equation:  BMW relations for E G Y 1}
\begin{split}
&YGY =  G + (q\inv - q) (1 - E),  \text{ and }\\
&G E = EG = \rho\inv E.
\end{split}
\end{equation}
   We give some details for the computation
of $YGY$;   applying $YGY$  to $\unitvector e_j$  gives
$$
Y G Y \, \unitvector e_j = 
(q\inv - q)  \sum_k  \frac{u_j u_k (\gamma_j - \delta(j, k))}{1 - u_j u_k} \, \unitvector e_k.
$$
If we write $u_j u_k = (u_j u_k -1)  + 1$  in the numerator and simplify, we get
\begin{equation*}
\begin{split}
 (q\inv - q)   &(- \gamma_j  \sum_k  \unitvector e_k +   \unitvector e_j)  + 
 (q\inv - q)  \sum_k \frac{ (\gamma_j - \delta(j, k))}{1 - u_j u_k} \, \unitvector e_k \\
 &=  (q\inv - q) (-E + 1) \, \unitvector e_j  + G   \, \unitvector e_j.
\end{split}
\end{equation*}
Taking into account  (\ref{equation:  BMW relations for E G Y 1}), we see that $YGY$ commutes with both
$E$ and $G$,
\begin{equation} \label{equation:  BMW relations for E G Y 3}
\begin{split}
&G Y G Y = Y G Y G, \text{ and } \\
&E Y G Y = Y G Y E.\\
\end{split}
\end{equation}
Similarly, we compute
$$
Y G Y E\, \unitvector e_j =  (q\inv - q)  \gamma_j \sum_k \sum_i  \frac{u_i u_k (\gamma_i - \delta(i, k))}{1 - u_i u_k} \, \unitvector e_k.
$$
Simplfying this as above, we get
$$
 (q\inv - q)  \gamma_j  (- \sum_i \gamma_i + 1) \sum_k \unitvector e_k +
 (q\inv - q)  \gamma_j \sum_k \sum_i  \frac{ (\gamma_i - \delta(i, k))}{1 - u_i u_k} \, \unitvector e_k.
$$
Using $\sum_i \gamma_i = \delta_0$  and   (\ref{equation:  equation for gammas00}), 
we can simplify this to
$$
((q\inv -  q)(- \delta_0 + 1)  + \rho\inv)\, \gamma_j   \sum_k \unitvector e_k = \rho \gamma_j  \sum_k \unitvector e_k = \rho E\, \unitvector e_j
$$
Thus, we have
\begin{equation}\label{equation:  BMW relations for E G Y 4}
E Y G Y = Y G Y E = \rho E.
\end{equation}

Set $m = \sum_i \unitvector e_i$.  We claim that $\{m, Y m, \dots, Y^{r-1} m\}$ is a basis of $M$.
In fact, if $\sum_{j = 0}^{r-1} \alpha_j Y^j m = 0$,  then  $\sum_{i = 1}^r (\sum_{j = 0}^{r-1}  \alpha_j u_i^j) \unitvector e_i = 0$.  Since the $u_i$ are distinct,  the Vandermonde matrix
$(u_i^j)$ is invertible,  so $\alpha_j = 0$ for all $j$.  This shows that $\{m, Y m, \dots, Y^{r-1} m\}$ is linearly independent, and therefore a basis of $M$.

We claim that  $G Y G Y$ is the identity transformation.  By the last paragraph,  it is enough to show that  $G Y G Y (Y^k m) = Y^k m$  for $0 \le k \le r-1$.  
Note that $YGYG E = E$, by  (\ref{equation:  BMW relations for E G Y 1}) and
(\ref{equation:  BMW relations for E G Y 4}), and also $E m = \delta_0 m$.  
According to 
(\ref{equation:  BMW relations for E G Y 3}),   $Y$ commutes with $GYG$, so 
\begin{equation*}
\begin{split}
G Y G Y (Y^k m) &=  Y^{k} YG Y G m = (1/\delta_0)  Y^{k} YG Y G  E m \\
&= (1/\delta_0)  Y^{k}  E m  =  Y^k m.
\end{split}
\end{equation*}
Using  (\ref{equation:  BMW relations for E G Y 1}),  we have 
\begin{equation} \label{equation:  BMW relations for E G Y 5}
G\inv = YGY =  G + (q\inv - q) (1 - E).
\end{equation}
It follows from  (\ref{equation:  BMW relations for E G Y 0})--(\ref{equation:  BMW relations for E G Y 5}) that $E, G, Y$ satisfy the relations of the cyclotomic BMW algebra with parameters
$\rho$,  $q$,  $\delta_j$ ($j \ge 0$), and $u_1, \dots, u_r$.  In particular, since we are working over a field, and $E \ne 0$,  the parameters must be weakly admissible;  see the statement just before Definition \ref{definition: weak admissibility}.
\end{proof}

\begin{definition} \label{definition: u--admissible}
 Let $S$ be an integral ground ring with parameters $\rho$, $q$,  $\delta_j$ ($j \ge 0$), and $u_1, \dots, u_r$.  Say that the parameters  are
{\em $u$-admissible}  or that $S$ is $u$--admissible if:
\begin{enumerate}
\item  $q - q\inv \ne 0$, 
 $u_i \ne u_j$ for $i \ne j$  and $u_i u_j \ne 1$  for
any $i, j$, and
\item 
the equivalent conditions of Theorem \ref{theorem: equivalent conditions for u--admissibility} hold.
\end{enumerate}
\end{definition}

\begin{corollary} \label{remark:  when is an admissible ring u--admissible}
Let $S$ be an integral ground ring with parameters $\rho$, $q$,  $\delta_j$ ($j \ge 0$), and $u_1, \dots, u_r$.  Then  $S$ is $u$--admissible if, and only if,  $S$ is admissible,  the $u_i$ are distinct and $u_i u_j \ne 1$  for all $i, j$.
\end{corollary}

\begin{proof} This follows from Theorems \ref{theorem: equivalent conditions for admissibility} and \ref{theorem: equivalent conditions for u--admissibility}. 
\end{proof}

\begin{remark} \label{remark:  changes of parameters for u--admissible rings}
Let $S$ be a ground ring with $u$--admissible parameters  $\rho$, $q$, $\delta_j$   ($j \ge 0$),  and
$u_1, \dots, u_r$.  Then  
$$
\rho, -q\inv, \delta_j \  (j \ge 0),  \text{ and }  u_1, \dots, u_r
$$
and
$$
-\rho, -q, \delta_j \  (j \ge 0),  \text{ and }  u_1, \dots, u_r
$$
are also sets of $u$--admissible parameters.   Moreover,  the quantities $\gamma_j$ computed from any of these sets of parameters according to  (\ref{equation: formula for gammas00}) are the same.
\end{remark}

\begin{remark}  \label{remark: new formulas for gammas}
 Let $S$ be an integral ground ring with $u$--admissible parameters $\rho$, $q$,  $\delta_j$ ($j \ge 0$), and $u_1, \dots, u_r$.  By  Remarks
\ref{remark:  solve for rho if S is an integral domain} and \ref{remark:  changes of parameters for u--admissible rings},  we can assume $\rho = \prod_j u_j$ if
$r$ is odd and $\rho = q\inv \prod_j u_j$  if $r$ is even.    Using this,  we can rewrite the formula for $\gamma_j$  (\ref{equation: formula for gammas00})  as
\begin{equation}\label{equation:  new formula for gamma odd case}
\gamma_j =  \frac{\rho (u_j -q\inv)(u_j + q)}{u_j^2(q - q\inv)} \prod_{\ell \ne j} \frac{u_j - u_\ell\inv}{u_j - u_\ell}  \quad \text{ if $r$ is odd},
\end{equation}
and
\begin{equation}\label{equation:  new formula for gamma even case}
\gamma_j =  \frac{\rho (u_j -q)(u_j + q)}{u_j^2(q - q\inv)} \prod_{\ell \ne j} \frac{u_j - u_\ell\inv}{u_j - u_\ell}  \quad \text{ if $r$ is even}.
\end{equation}
Note that if $u_j \not\in \{\pm q\powerpm\}$,  then $\gamma_j \ne 0$.
\end{remark}

\subsection{Some generating functions} \label{subsection:  symmetric functions}
We will investigate some generating functions that appear naturally in the study of cyclotomic BMW algebras.
Let $u_1,  \dots, u_r$,  $q$,  and $\rho$  be indeterminants over $\Z$.  We will need to consider rational functions $\gamma_j$  in these variables that are solutions to the system of linear equations:
\begin{equation*}
\sum_j  \frac{1}{1-u_i u_j}\,  \gamma_j =    \frac{1}{1-u_i^2}  +  \frac{1}{\rho (q\inv - q)} \quad (1 \le i \le r).
\end{equation*}

\begin{lemma}  \label{lemma:  gammas from system of equations}
Let $F$ be field and let $u_1, \dots, u_r$  be distinct non--zero elements of $F$ with $1 \ne u_i u_j$  for all $i, j$.
\begin{enumerate}
\item  The matrix  $$\left (\frac{1}{1 - u_i u_j} \right )_{1 \le i, j \le r}$$ is invertible.
\item  The unique solution to the linear system of equations
$$  \sum_j   \frac{1}{1 - u_i u_j}   \gamma_j^{(1)}  =   1   \quad (1 \le i \le r)$$
is
$$
\gamma_j^{(1)} =  (1 - u_j^2) \prod_{\ell \ne j}  u_\ell \prod_{\ell \ne j}  \frac{(u_\ell u_j -1)}{u_j - u_\ell}.
$$
\item  The unique solution to the linear system of equations
$$
\sum_j  \frac{1}{1 - u_i u_j}   \gamma_j^{(2)} = \frac{1}{1- u_i^2}  \quad (1 \le i \le r)
$$
is
$$
 \gamma_j^{(2)} = \prod_{\ell \ne j}  \frac{u_\ell u_j -1}{u_j - u_\ell}  \quad \text{ if $r$ is odd},
$$
and 
$$
 \gamma_j^{(2)} = - u_j  \prod_{\ell \ne j}  \frac{u_\ell u_j -1}{u_j - u_\ell} \quad   \text{ if $r$ is even}.
$$
\end{enumerate}
\end{lemma}

\begin{proof}  Recall  Cauchy's determinant identity ~\cite{Kratt}, Section 2.1, 
$$
\det\left( \frac{1}{1- x_i y_j} \right ) =  \frac{\displaystyle \prod_{i < j} (x_i - x_j) (y_i - y_j)}  {\displaystyle \prod_{i, j} (1 - x_i y_j)},
$$
which specializes to 
$$
\det\left( \frac{1}{1- u_i u_j} \right ) =  \frac{\displaystyle \prod_{i < j} (u_i - u_j)^2}  {\displaystyle \prod_{i, j} (1 - u_i u_j)}.
$$
Since $u_i \ne u_j$ for $i \ne j$,  this implies statement (1).

For part (2),  we have to prove the identity, for all $i$  ($1 \le i \le r$):
$$
\sum_j    \frac{1}{1 - u_i u_j}   (1 - u_j^2) \prod_{\ell \ne j}  u_\ell \prod_{\ell \ne j}  \frac{(u_\ell u_j -1)}{u_j - u_\ell} = 1,  
$$
or
\begin{equation}
\label{equation: identity determining first part of gammas}
\sum_j \  \prod_{\ell \ne i}( u_\ell u_j -1)  \prod_{\ell \ne j}  \frac{u_\ell}{u_j - u_\ell}  = 1.
\end{equation}
By the principal of permanence of identities (see ~\cite{artin}, p. 456)  it suffices to prove this when
the $u_i$  are distinct non--zero complex numbers satisfying $u_i u_j \ne 1$  for all $i, j$.
Put
$$
f_i(\zeta) =  (1/\zeta) \prod_{\ell \ne i} (\zeta u_\ell -1)  \prod_{\ell = 1}^r \frac{u_\ell}{\zeta - u_\ell}.
$$
Check that the residue of $f_i(\zeta)$  at $\zeta = u_j$  is the summand on the left side of 
 (\ref{equation: identity determining first part of gammas}).  The other possible poles of $f_i(\zeta)$ are at $\zeta = 0$ and at $\zeta = \infty$.    The residue at $\zeta = 0$ is $-1$  and the residue at $\zeta = \infty$ is zero.  The identity of  (\ref{equation: identity determining first part of gammas}) follows, since the sum of the residues of 
$f_i(\zeta)$  at its finite poles and at $\infty$  is zero.

For part (3), we have to prove the identity for all $i$  ($1 \le i \le r$):
\begin{equation}
\label{equation:  identity determining second part of gammas,  odd case}
\sum_j  \frac{1 - u_i^2}{1- u_i u_j}  \prod_{\ell \ne j}  \frac{u_\ell u_j -1}{u_j - u_\ell} = 1  \quad \text{ if $r$ is odd},
\end{equation}
and 
\begin{equation}
\label{equation:  identity determining second part of gammas,  even case}
\sum_j  \frac{1 - u_i^2}{1- u_i u_j} (-u_j) \prod_{\ell \ne j}  \frac{u_\ell u_j -1}{u_j - u_\ell} = 1  \quad \text{ if $r$ is even}.
\end{equation}
Again, we may assume that the $u_i$ are distinct non--zero complex numbers.  When $r$ is odd, we put
$$
f_i(\zeta) = \frac{1- u_i^2}{(\zeta^2 - 1)(1 - \zeta u_i)} \prod_{\ell = 1}^r \frac{\zeta u_\ell - 1}{\zeta - u_\ell},
$$
and observe that the residue of $f_i(\zeta)$  at $\zeta = u_j$ is the summand on the left side of 
 (\ref{equation:  identity determining second part of gammas,  odd case}).  The other possible poles of
$f_i(\zeta)$  are at $\zeta = \pm 1$ and $\zeta = \infty$.   The residue at $\zeta = 1$ is $(-1/2)(1 + u_i)$,  the
residue at $\zeta = -1$ is $(-1/2)(1 - u_i)$,   and the residue at $\zeta = \infty$ is zero.  Since the sum of these residues   is $-1$,  Equation (\ref{equation:  identity determining second part of gammas,  odd case}) follows.
The proof of  (\ref{equation:  identity determining second part of gammas,  even case}) is similar, and we omit it.
\end{proof}

\begin{remark}  \label{proof of formula for gammas}
Lemma \ref{corollary:  gammas defined by linear equations} from Section \ref{subsection: u--admissibility} is a corollary  of Lemma \ref{lemma:  gammas from system of equations}.  Namely,  the unique solution to the system of equations
\begin{equation*}
\sum_j  \frac{1}{1-u_i u_j}\,  \gamma_j =    \frac{1}{1-u_i^2}  +  \frac{1}{\rho (q\inv - q)} \quad (1 \le i \le r)
\end{equation*}
is
$$
\gamma_j =  \frac{1}{\rho (q\inv -q)} \gamma_j^{(1)}  +  \gamma_j^{(2)}.
$$
\end{remark}

Fix a natural number $r$.  Let $u_1,  \dots, u_r$ and $t$  be indeterminants over $\Z$,
and define
\begin{equation}\label{equation: definition of G(t)}
G(t) = G(u_1, \dots, u_r; t) = \prod_{\ell = 1}^r  \frac{t - u_\ell}{t u_\ell -1}.
\end{equation}
Let $\mu_a = \mu_a(u_1, \dots, u_r)$  denote the coefficients of the formal power series expansion of $G(t)$,
\begin{equation} \label{equation:  mu's as coefficients of G(t)}
G(t) =  \sum_{a = 0}^\infty  \mu_a  t^a.
\end{equation}
Notice that each $\mu_a$  is a symmetric polynomial in $u_1,  \dots, u_r$  and that
$G(t\inv)  =  G(t)\inv$.

Now suppose that $\rho$ and $q$ are additional indeterminants, and for $1 \le j \le r$,  define rational functions $\gamma_j$  by  (\ref{equation: formula for gammas00}).  Moreover, for $a \ge 0$,  define $
\delta_a =  \sum_{j = 1}^r  \gamma_j  u_j^a
$.  
Let $Z_1(t)$  be the generating function for the $\delta_a$,
\begin{equation} \label{equation: definition of Z1}
\begin{aligned}
Z_1(t) & =  \sum_{a \ge 0}  \delta_a t^{-a} =  \sum_{a \ge 0} \sum_{j=1}^r \gamma_j u_j^a t^{-a} \\
&=  \sum_{j=1}^r    \gamma_j \sum_{a \ge 0} u_j^a t^{-a} =   \sum_{j=1}^r    \gamma_j  \frac{t}{t - u_j}.
\end{aligned}
\end{equation}

In the following, we use the notation  $$\delta(P) = \begin{cases}
1 &\text{if $P$  is true }\\
0 &\text{if $P$ is false}.
\end{cases}
$$

\begin{lemma}  \label{lemma:  calculation of Z1}   Let $p$  denote $\prod_{j=1}^r  u_j$.
\begin{enumerate}
\item  \quad  $\displaystyle Z_1(t) =  \frac{1}{\rho(q\inv -q)}   +  \frac{t^2}{t^2 -1} + A(t) \  G(t\inv)$,
where
$$ 
A(t) =
\begin{cases}
\displaystyle
{\rho\inv p }/{(q - q\inv)}\ +\  {t}/{(t^2-1)}  \quad &\text{if $r$ is  odd, and } \\
\displaystyle
{\rho\inv p}/{(q - q\inv)} \ -\  {t^2}/{(t^2-1)} \quad &\text{if $r$ is  even}.
\end{cases}
$$

\item

 If $r$ is odd, then
for $a \ge 0$,  
$$
\begin{aligned}
\delta_a =  \delta(a = 0)\ & {\rho\inv  }/{(q - q\inv)}
+ \delta(\text{$a$ is even})  \\ &+     \mu_a \,{\rho\inv p }/{(q - q\inv)}  + \mu_{a-1}  + \mu_{a-3}  + \cdots  .
\end{aligned}
$$

\item
If $r$ is even,  then
for $a \ge 0$,  
$$
\begin{aligned}
\delta_a =  \delta(a = 0)\ & {\rho\inv  }/{(q - q\inv)}
+ \delta(\text{$a$ is even})  \\ &+     \mu_a  \,{\rho\inv p }/{(q - q\inv)} - \mu_a   - \mu_{a-2}  - \mu_{a-4}  + \cdots  .
\end{aligned}
$$

\item  $\displaystyle \delta_0 =  \frac{1-  p^2 }{\rho (q\inv -q)}  
+ 1  -
\begin{cases}
0  &\text{ if $r$ is odd} \\
p  &\text{ if $r$ is even}.
\end{cases}
 $
 \item  For all $a \ge 0$,   $\delta_a$ is an element of the ring
 $$\Z[u_1,  \dots, u_r,  \rho\inv,  (q\inv -q)\inv],$$
 and is symmetric in $u_1, \dots, u_r$.
\end{enumerate}
\end{lemma}

\begin{proof}
To prove the identity of part (1),  it suffices to suppose that the quantities $u_1,  \dots, u_r$, $t$, $q$,  and $\rho$  are complex numbers,  algebraically independent over $\Q$.  
Set $$
h(\zeta) = G(\zeta\inv)  \left(
\zeta\inv \frac{p}{  \rho (q - q\inv)}
+
\begin{cases}
{1}/{(\zeta^2 -1)}  &\text{if $r$  is odd} \\
{-\zeta}/{(\zeta^2 -1)}  &\text{if $r$  is even}
\end{cases}
\right ).
$$
It is straighforward to check that the residue of $h(\zeta) t/ (t-\zeta)$  at $\zeta = u_j$
is $$\gamma_j   t/ (t-u_j).$$
Thus,  $Z_1(t)$  is the sum of the residues of $h(\zeta) t/ (t-\zeta)$  at $u_1,  \dots, u_r$.
The other possible poles of $h(\zeta) t/ (t-\zeta)$  are at $0$,   $\pm 1$,   $\infty$,  and $t$.
Independent of the parity of $r$,  the residue of  $h(\zeta) t/ (t-\zeta)$ at $\zeta = 0$ is
$
1/(\rho(q - q\inv)),
$
the residue at $\zeta = 1$ is
$
-(1/2)(t/(t-1)),
$
the residue at $\zeta = -1$ is
$
-(1/2)(t/(t+1)),
$
and the residue at $\infty$ is zero.  For $r$ odd,  the residue of $h(\zeta) t/ (t-\zeta)$ at $\zeta = t$ is
$$
G(t\inv) \left(
p/ (\rho(q\inv -q))   -  t/(t^2 -1)
\right ).
$$
For  $r$ even,  the residue of $h(\zeta) t/ (t-\zeta)$ at $\zeta = t$ is
$$
G(t\inv) \left(
p/ (\rho(q\inv -q))   +  t^2/(t^2 -1)
\right ).
$$
$Z_1(t)$  is the negative of the sum of the residues of $h(\zeta) t/ (t-\zeta)$ at $0$, $\pm 1$, 
$\infty$,  and $t$.   This yields the formula of part (1).

For part (2), write
$$
\begin{aligned}
Z_1(t\inv) &= \sum_{a \ge 0}  \delta_a  t^a 
= 
 \frac{1}{\rho(q\inv -q)}   +  \frac{1}{1-t^2 }  \\
&+  G(t) \left( \frac{p}{\rho(q - q\inv)} +
\begin{cases}  \displaystyle
 {t}/({1-t^2})  &\text{if $r$ is  odd}\\
\displaystyle {-1}/({1-t^2})  &\text{if $r$ is  even}
\end{cases}
\right ).
\end{aligned}
$$
One can now read off the coefficient of $t^a$ on the right hand side.

Part (3) follows from part (2),  using that $\mu_0 = G(0) = p$,   and part (4) follows from part (2)  as well.
\end{proof}

Note that if $F$ is a field with $u$--admissible parameters $\rho$, $q$, $\delta_j$,  and
$u_1, \dots, u_r$,  and  $Z_1(t) = \sum_{a \ge 0}  \delta_a t^{-a}$,  then  Lemma 
\ref{lemma:  calculation of Z1} applies,  but we also can assume that $\rho = p$ if $r$ is odd,
and $\rho = q\inv p$ if $r$ is even.  So the formulas of the lemma simplify as follows:

\begin{equation} \label{equation:  formula for Z1 in u admissible case}
\begin{aligned}
\displaystyle Z_1(t) &=  \frac{1}{\rho(q\inv -q)}   +  \frac{t^2}{t^2 -1} + A(t) \  G(t\inv), \quad 
 \text{where}\\
A(t) &=
\begin{cases}  \displaystyle{1}/{(q - q\inv)}\ +\  {t}/{t^2-1}  \quad &\text{if $r$ is  odd, and } \\
 \displaystyle
{q}/{(q - q\inv)} \ -\  {t^2}/{t^2-1} \quad &\text{if $r$ is  even.}
\end{cases}
 \\
\end{aligned}
\end{equation}

\subsection{The universal admissible ring}

We begin by constructing a universal $u$--admissible integral domain.  Let $\qbold$, $\ubold_1$, \dots, $\ubold_r$   be algebraically independent indeterminants over $\Z$.    Let $S_0$ be the Laurent polynomial ring in $\qbold$, $\ubold_1$, \dots, $\ubold_r$,  with
$ (\qbold\inv - \qbold)\inv$ adjoined,
$$
S_0 = \Z[\qbold^{\pm 1}, \ubold_1^{\pm 1}, \dots, \ubold_r^{\pm 1}][  (\qbold\inv - \qbold)\inv].
$$
Let $\bm p = \prod_j \ubold_j$. 
If $r$ is odd,  put $\rhobold = \bm p$,  and if $r$ is even,  put  $\rhobold = \qbold\inv \bm p$.
Define $\bm \gamma_j$  for $1 \le j \le r$ (in the field of fractions of $S_0$) by  (\ref{equation: formula for gammas00}), and set 
$\deltabold_a = \sum_{j = 1}^r \bm \gamma_j \ubold_j^a$.    
According to Lemma \ref{lemma:  calculation of Z1},   $\deltabold_a \in S_0$,  and, moreover, 
\begin{equation*} 
\deltabold_0 = 1 + \frac{1 - \bm p^2} {\rhobold(\qbold\inv - \qbold)} -
\begin{cases}
0 &\text{if $r$ is odd} \\
\bm p  &\text{if $r$  is even}.
\end{cases} 
\end{equation*}
  Taking into account that $\rhobold = \bm p$  or
$\rhobold = \qbold\inv \bm p$,  we have $\deltabold_0 \ne 0$  by the algebraic independence of $\qbold$  and $\ubold_1,  \dots, \ubold_r$.    Moreover,  one can check that \break $\rhobold\inv -  \rhobold = (\qbold\inv - \qbold)(\deltabold_0 - 1)$.
Finally,   define
$
\overline S = S_0[\deltabold_0\inv].
$
Then  $\overline S$ is $u$--admissible.  Note that the field of fractions of $\overline S$ is $\Q(\qbold, \ubold_1, \dots, \ubold_r)$.

Let   $S$  be another $u$--admissible integral domain with parameters $\rho$,  $q$, etc.
Using Remark \ref{remark:  changes of parameters for u--admissible rings}, we can assume that  $\rho = \prod_j u_j$ if
$r$ is odd and $\rho = q\inv \prod_j u_j$  if $r$ is even.  
Because of the algebraic independence of $\qbold$, $\ubold_1$, \dots, $\ubold_r$, there is a ring homomorphism $\varphi : \overline S \rightarrow S[(q - {q}\inv)\inv]$
mapping $\qbold \mapsto q$  and $\ubold_i \mapsto u_i$  for each $i$.  It follows that $\varphi: \rhobold \mapsto \rho$.  Since $\overline S$ and $S$  are $u$--admissible,  we have 
$\deltabold_a \mapsto \delta_a$  for $a \ge 0$  as well. 

We have shown the following result:

\begin{lemma}  There is a ``universal" integral domain $\overline S$ with $u$--admissible parameters $\rhobold$, $\qbold$,  $\deltabold_j$, ($j \ge 0$),  $\ubold_1, \dots, \ubold_r$, with the properties:
\begin{enumerate}
\item  $\qbold$, $\ubold_1$, \dots, $\ubold_r$ are algebraically independent over $\Z$.
\item $\overline S$  is generated as a unital ring by 
 $\qbold^{\pm 1}$,  $(\qbold\inv - \qbold)\inv$,  $\ubold_1^{\pm 1}, \dots, \ubold_r^{\pm 1}$, and $\deltabold_0\inv$.
 \item The field of fractions of $\overline S$ is $\Q(\qbold, \ubold_1, \dots, \ubold_r)$.
 \item  Whenever $S$ is an integral domain with $u$--admissible
parameters $\rho$, $q$, etc.,    there exists a  morphism of ground rings (see Definition \ref{definition: parameter preserving})
$$\phi : \overline S \rightarrow S[(q-q\inv)\inv].$$
\end{enumerate}
\end{lemma}

Next we consider the construction of  universal admissible rings.  Let 
$\qbold$, $\rhobold$, $\deltabold_0, \dots, \deltabold_{r-1}$,  $\ubold_1, \dots, 
\ubold_r$  be indeterminants over $\Z$. 
 Let $\abold_j$  be the signed elementary symmetric polynomials in the variables $\ubold_i$.
For $a \ge r$, define $\deltabold_a$ by the recursion  $\deltabold_a = 
 -\sum_{j=0}^{r-1} \abold_j \deltabold_{a-r+j}$. 
 
 Set 
$$
\A = \Z[\qbold\powerpm, \rhobold\powerpm, \deltabold_0\powerpm, \deltabold_1, \dots, \deltabold_{r-1},  \ubold_1\powerpm, \dots, \ubold_r\powerpm].
$$
If $r$ is odd, set $f_+ = \rhobold + \abold_0$  and $f_- = \rhobold - \abold_0$.  If $r$ is even, set  $f_+ = \rhobold - \qbold\inv \abold_0$  and $f_- = \rhobold + \qbold \abold_0$. 
For $r$ of either parity,  let $f_0 = f_+ f_-$.  For $\eps \in \{0, +, -\}$,  let 
$J_\epsilon$  be the ideal in $\A$ generated by $f_\eps$, $\rhobold\inv -  \rhobold - (\qbold\inv - \qbold)(\deltabold_0 - 1)$, and the elements
\begin{equation*} 
\begin{aligned}
&\rhobold(\abold_\ell - \abold_{r-\ell}/\abold_0) \ + \\&(\qbold-\qbold\inv)\bigg [ \sum_{j = 1}^{r - \ell} \abold_{j+\ell} \deltabold_j 
-  \sum_{j = \max(\ell + 1, \lceil r/2 \rceil)}^{\lfloor (\ell + r)/2 \rfloor} \abold_{2j - \ell} 
+  \sum_{j =  \lceil \ell/2 \rceil}^{\min(\ell, \lceil r/2 \rceil -1)} \abold_{2j - \ell} \bigg ],
\end{aligned}
\end{equation*}
for $1 \le \ell \le r-1$.
Note that $J_0 \subseteq  J_+ \cap J_-$.     Set $R_\eps = \A/J_\eps$.
Since $J_0 \subseteq J_\pm$,  we have maps $R_0 \to R_\pm$  defined by
$x + J_0 \mapsto x + J_\pm$.  Moreover,  we can define an automorphism  $\theta$ of $\A$ taking
$J_+$ to $J_-$.   (For $r$ odd,   require $\theta(\rhobold) = -\rhobold$,  $\theta(\qbold) = -\qbold$,  $\theta(\deltabold_j) = \deltabold_j$  and $\theta(\ubold_i) = \ubold_i$.  For $r$ even,  require $\theta(\rhobold) = \rhobold$,  $\theta(\qbold) = -\qbold\inv$,  $\theta(\deltabold_j) = \deltabold_j$  and $\theta(\ubold_i) = \ubold_i$.)   Then $\theta$ induces an isomorphism $\theta: R_+ \to R_-$.

The rings $R_\eps$  have the properties of an admissible ring,  except that  the image of $(\qbold - \qbold\inv)$ in $R_\eps$   might be a zero--divisor.   Consequently, we will introduce quotients of $R_\eps$ in which the image
of $(\qbold - \qbold\inv)$ is not a zero divisor.

If $S$ is an admissible integral domain,  then there is a  morphism of ground rings from $R_0$ to $S$, and this homomorphism must  factor through $R_+$ or $R_-$ according to Remark \ref{remark:  solve for rho if S is an integral domain}.  In particular,  there is a morphism of ground rings from 
$R_0$ to the universal $u$--admissible ring $\overline S$ that factors through $R_+$.   It follows from this that
$(\qbold - \qbold\inv)^n + J_0 \ne 0$ in $R_0$  and   $(\qbold - \qbold\inv)^n + J_+ \ne 0$ in $R_+$, for all natural numbers $n$.     Using the isomorphism $\theta: R_+ \to R_-$,  we have also that
$(\qbold - \qbold\inv)^n + J_- \ne 0$ in $R_-$ for all natural numbers $n$.

For $\eps \in \{0, +, -\}$, let
$$
P_\eps = \{ x + J_\eps \in R_\eps :  \exists n \ge 1 \text{ such that }  (x + J_\eps) ( \qbold - \qbold\inv + J_\eps)^n = 0 \}.
$$
Then  $P_\eps$ is an ideal in $R_\eps$.   Set $\overline R_\eps = R_\eps/P_\eps$.  Let $\bar\rho$,  $\bar q$, etc. denote the images of the parameters in $\overline R_\eps$.
Since the quotient maps $R_0 \to R_\pm$ take  $P_0$ into $P_\pm$,  they induce maps $\overline R_0  \to \overline R_\pm$.  Moreover,  the isomorphism $\theta: R_+ \to R_-$  maps $P_+$ to $P_-$,  so induces an isomorphism $\theta: \overline R_+ \to \overline R_-$.

If $S$ is  any admissible ring and  $\varphi : R_\eps \rightarrow S$  is a
morphism of ground rings from $R_\eps$ to $S$, then $\varphi(P_\eps) = 0$.  Hence $\varphi$ induces a morphism of ground rings from 
$\overline R_\eps$ to $S$.    In particular,  we have a morphism of ground rings from $R_0$ to the universal $u$--admissible ring $\overline S$,  that factors through $R_+$,  and this map induces a morphism of ground rings from $\overline R_0$ to $\overline S$,  factoring through $\overline R_+$.
It follows from this that $\bar q - \bar q\inv \ne 0$ in $\overline R_\eps$  ($\eps = 0, +$).  
Using the isomorphism $\theta: \overline R_+ \to \overline R_-$,  we see that
$\bar q - \bar q\inv \ne 0$ in $\overline R_-$  as well.

Finally we check that $\bar q - \bar q\inv$ is not a zero--divisor in $\overline R_\eps$, for 
$\eps \in \{0, +, -\}$.   Let $\bar x \in \overline R$ such that $\bar x (\bar q - \bar q\inv) = 0$.   Let $x$ be a preimage of $\bar x$ in $R_\eps$.   Then
$x (\qbold  - \qbold\inv + J_\eps) \in P_\eps$.  Hence there exists an $n \ge 1$ such that
$$x ((\qbold  - \qbold\inv) + J_\eps) ((\qbold  - \qbold\inv) + J_\eps)^n = 0.$$  But this means that
$x \in P_\eps$, so $\bar x = 0$.

We claim that $\overline R_+$ is an integral domain.

Since there is a morphism of ground rings from $\overline R_+$ to the universal $u$--admissible ring $\overline S$,  it follows that the parameters $\bar q$,  $\bar u_1$, \dots, $\bar u_r$  in $\overline R_+$  are algebraically independent over $\Z$.  
Consider $\overline R_+[(\bar q - \bar q\inv)\inv]$  and its subring
$$A_0 = \Z[\bar q\powerpm, \bar u_1\powerpm, \dots, \bar u_r\powerpm][ (\bar q - \bar q\inv)\inv].$$  
$A_0$ is an integral domain,  since $\bar q$,  $\bar u_1$, \dots, $\bar u_r$  are algebraically independent over $\Z$.  According to Corollary \ref{corollary:  rho and deltas contained in subring},  $A_0$ also contains
$\bar\rho\powerpm$ and $\bar \delta_0, \bar \delta_1, \dots, \bar \delta_{r-1}$.
Therefore,  $$A_0[\bar \delta_0\inv] = \overline R_+[(\bar q - \bar q\inv)\inv].$$
  Thus,
$\overline R_+[(\bar q - \bar q\inv)\inv]$ and its subring $\overline R_+$ are integral domains.
According to Corollary \ref{remark:  when is an admissible ring u--admissible},  $\overline R_+$ is $u$--admissible.  It is now clear that  $\overline R_+[(\bar q - \bar q\inv)\inv]$ is isomorphic to the universal $u$--admissible ring $\overline S$  by a isomorphism of ground rings.

We summarize this discussion with the following theorem:

\begin{theorem}\mbox{} \label{theorem:  universal admissible rings}  
There exist a universal admissible ring $\overline R_0$  and a universal admissible integral domain $\overline R_+$  with the following properties:
\begin{enumerate}
\item  There is  morphism of ground rings from $\overline R_0$ to $\overline R_+$.
\item   $\overline R_0$   (and hence also $\overline R_+$) are generated as  unital rings by the parameters
 $\bar q^{\pm 1}$, $\bar \rho \powerpm$, $\bar \delta_0\powerpm$, $\bar \delta_1$, \dots $\bar \delta_{r-1} $,$\bar u_1^{\pm 1}, \dots, \bar u_r^{\pm 1}$.
 \item The parameters $\bar q, \bar u_1, \dots, \bar u_r$ of  
$\overline R_+$  (hence also of  $\overline R_0$)  are algebraically independent over $\Z$.
\item Whenever $S$ is a ring with admissible
parameters,    there exists a  morphism of ground rings from $ \overline R_0$ to  $S$.  If $S$ is also an integral domain, then the morphism of ground rings can be chosen to factor through $\overline R_+$.
\item $\overline R_+[(\bar q - \bar q\inv)\inv]$ is isomorphic to the universal $u$--admissible ring, by an  isomorphism of ground rings.
\end{enumerate}
\end{theorem}

\subsection{More generating functions}  \label{subsection:  more generating functions}

 Let $F$ be a field with $u$--admissible parameters
$\rho$, $q$,  $\delta_j$,  $u_1, \dots, u_r$.   Let $\kt {n, r}$  denote
$\kt {n, F,r }(u_1, \dots, u_r)$.  Recall that $\eps_n$ denotes the conditional expectation from
$\kt {n, r}$  to $\kt {n-1, r}$, and $\eps$ denotes the canonical trace on $\kt {n, r}$.    For $n \ge 1$, write
$Y_n$  for $\varphi(y_n)$.  Note that $\eps_n(Y_n^a) \in Z(\kt {n-1, r})$ for $a \ge 0$,   since
$Y_n^a$ commutes with $\kt {n-1, r}$.  Define $\omega_n^{(a)}$ by
$\omega_n^{(a)} = \delta_0 \,\eps_n(Y_n^a)$, so 
$$
E_n Y_n^a E_n = \omega_n^{(a)} E_n.
$$
Define
\begin{equation}  \label{equation: definition of Zk}
Z_n(t) = \sum_{a \ge 0}   \omega_n^{(a)} t^{-a},
\end{equation}
a formal power series with coefficients in the center of $\kt {n-1, r}$.  This is consistent with the previous definition of $Z_1$ in  \ref {equation: definition of Z1}.   Moreover,  define
\begin{equation}  \label{equation: definition of Qk}
Q_n(t) =  Z_n(t) - \rho\inv/(q\inv - q)  -  t^2/(t^2 -1).
\end{equation}
In particular
\begin{equation}  \label{equation: expression for Q1 of t}
Q_1(t) =  G(t\inv)  A(t),
\end{equation}
according to  (\ref{equation:  formula for Z1 in u admissible case}).

We have the following remarkable recursion for $Q_n(t)$,  which is due to 
Beliakova  and  Blanchet  ~\cite{blanchet-beliakova},  Lemma 7.4.

\begin{proposition} \label{proposition:  recursion for Q's}
$$Q_{n+1}(t) = Q_n(t) \left (  \frac{(t-Y_n)^2 (t - q^{-2} Y_n\inv) (t- q^2 Y_n\inv)}
{(t-Y_n\inv)^2 (t - q^{-2} Y_n) (t- q^2 Y_n)}      \right ).$$
\end{proposition}

\begin{proof}  Beliakova and Blanchet prove this for the ordinary BMW  (or Kauffman tangle) algebras, that is,  for $r = 1$.  However,  their proof remains valid in the general case.
\end{proof}

\section{The generic structure of the cyclotomic BMW algebra}

In this section we obtain the ``generic structure"  of the cyclotomic BMW algebras over a field, following the method of Wenzl from ~\cite{Birman-Wenzl},  ~\cite{Wenzl-Brauer}, and  ~\cite{Wenzl-BCD}.  The argument also relies  on ideas from
Beliakova and Blanchet ~\cite{blanchet-beliakova}.

 Let $F$ be the field of fractions of the universal admissible integral domain $\overline R_+$.      We will show that  for all $n$  the cyclotomic BMW algebra $\bmw{n, F, r}$ over $F$ is isomorphic to the cyclotomic Kauffman tangle algebra $\akt {n, F, r}$ over $F$,  and  is split semisimple of dimension $r^n (2n-1)!!$.  The irreducible representations of $\bmw{n, F, r}$  are labelled
by $r$--tuples of Young diagrams of total size $n$, $n-2$, $n-4$,  \dots,  and the
dimension of the irreducible representation labelled by such an $r$--tuple $\mu$ is the number of up--down tableaux of length $n$ and shape $\mu$.   

\subsection{Semisimplicity and structure of $\bmw {n, F,r}$}

Let $S$  be a ground  ring.    For each $n$,  let $I_n$  be the two sided ideal in $\bmw{n, S, r}$   generated by $e_{n-1}$.   
Because of the relations  $e_j e_{j\pm1} e_j = e_j$,   the ideal $I_n$  is generated by any $e_i$  ($1 \le i \le n-1$)  or by all of them.   

\begin{lemma} \label{lemma:  I n plus 1 lemma}
For all $n \ge 1$,   $I_{n+1}  =  \bmw{n, S, r}  \,e_n \, \bmw{n, S, r}$.
\end{lemma}

\begin{proof}  This follows from the corresponding result for the affine BMW algebras, ~\cite{GH1},  Proposition 3.20.
\end{proof}

Let us recall the definition of the affine and cyclotomic Hecke algebras,   see ~\cite{ariki-book}.

\begin{definition}
Let $S$ be a commutative unital ring with an invertible element $q$.  The {\em affine Hecke algebra} 
$\ahec{n,S}(q^2)$ 
over $S$
is the $S$--algebra with generators $t_1,  g_1,  \dots, g_{n-1}$, with relations:
\begin{enumerate}
\item  The generators $g_i$ are invertible,  satisfy the braid relations,  and \break $g_i - g_i\inv =  (q - q\inv)$.
\item  The generator $t_1$ is invertible,  $t_1 g_1 t_1 g_1 = g_1 t_1 g_1  t_1$  and $t_1$ commutes with $g_j$  for $j \ge 2$.
\end{enumerate}
Let $u_1,  \dots, u_r$  be additional  elements in $S$.   The {\em cyclotomic Hecke algebra}
\break $\hec{n, S, r}(q^2; u_1,  \dots, u_r)$  is the quotient of the affine Hecke algebra $\ahec{n, S}(q^2)$ by the polynomial relation    $(t_1 - u_1) \cdots (t_1 - u_r) = 0$.
\end{definition}

\begin{remark} Since the generator $t_1$ can be rescaled by an arbitrary invertible element of $S$, only the ratios of the parameters $u_i$ have invariant significance in the definition of the cyclotomic Hecke algebra.
\end{remark}

\begin{lemma} \label{lemma:  hecke quotient of cyclotomic bmw}
 Let $S$ be a ground ring with parameters $\rho$,   $q$,  $\delta_j$  $(j \ge 0$), and $u_1,  \dots, u_r$.    For all $n \ge 1$,  the quotient of  the cyclotomic BMW algebra \break $\bmw{n, S, r}(u_1, \dots, u_r)$  by the ideal
$I_n$  is isomorphic to the cyclotomic Hecke algebra $\hec{n, S, r}(q^2; u_1, \dots, u_r)$.
\end{lemma}

\begin{proof}  Evident.
\end{proof}

\begin{lemma}  \label{lemma: isomorphism and semisimplicity for one strand}
Let $F$ be a field with $u$--admissible parameters $\rho$,   $q$,  $\delta_j$  $(j \ge 0$), and $u_1,  \dots, u_r$.    Assume that
$u_j \not\in \{\pm \, q\powerpm \}$ for all $j$.
\begin{enumerate}
\item  $\bmw {1, F, r} \cong \kt {1, F, r} \cong F^r$.
\item  Let $\displaystyle P_j = \prod_{\ell \ne j} \frac{Y_1 - u_\ell}{u_j - u_\ell}\in \kt {1, F, r}.$   Then
$E P_j E = \gamma_j E$,  and $\eps(P_j) = \gamma_j/\delta_0$, where $\gamma_j$ is defined by  (\ref{equation: formula for gammas00}).
\end{enumerate}
\end{lemma}

\begin{proof}
Let $\varphi: \bmw {2, F, r} \rightarrow \kt {2, F, r}$  be the surjective algebra homomorphism determined by
$\varphi(y_1) = Y_1$, $\varphi(e) = E$, $\varphi(g) = G$.   Define $\displaystyle  p_j =  \prod_{\ell \ne j} \frac{y_1 - u_\ell}{u_j - u_\ell}\in \bmw {1, F, r}$.  It is shown in the proof of Theorem \ref{theorem: equivalent conditions for u--admissibility} that
$e p_j e = \gamma_j e$.   Hence,  $E P_j E = \varphi(e p_j e) = \gamma_j E$.

By Remark  \ref{remark: new formulas for gammas},  $\gamma_j \ne 0$  for all $j$,  and by Lemma \ref{lemma:  E and G not zero},  $E \ne 0$.   Therefore,  $P_j \ne 0$.   The elements $P_j$  are mutually orthogonal non--zero idempotents,  so they are linearly independent.  Thus the restriction of $\varphi$ to $\bmw {1, F, r}$  is  an isomorphism. Moreover, 
$\kt {1, F, r} \cong F P_1 \oplus \cdots \oplus F P_r \cong F^r$,  as algebras.
\end{proof}

Let $\labold = (\la^{(1)}, \dots, \la^{(r)})$ be an $r$--tuple of Young diagrams.  The total size of $\labold$ is
$|\labold| = \sum_i |\la^{(i)}|$.    If $\mubold$ and $\labold$ are  $r$--tuples of Young diagrams of total size $f-1$  and $f$ respectively,  we write $\mubold \subset \labold$ if $\mubold$ is obtained from $\labold$ by removing one box from one component of $\labold$.

\begin{theorem}[\cite{ariki-book}]   \label{theorem:  cyclotomic hecke split semisimple}
Let $F$ be a field.   The cyclotomic Hecke algebra \break $\hec{n, F, r}(q; u_1, \dots, u_r)$  is split semisimple for all $n$  as long  as $q$ is not a proper root of unity and $u_i/u_j$  is not a power of $q$  for all
$i \ne j$.  In this case,   the simple components of $\hec{n, F, r}(q; u_1, \dots, u_r)$   are labeled by $r$--tuples of Young diagrams of total size $n$,   and a simple $\hec{n, F, r}$  module $V_\labold$ decomposes as a
$\hec{n-1, F, r}$  module as the direct sum of all $V_\mubold$  with $\mubold \subset \labold$.
\end{theorem}

In the following,  let $\Gamma_n$  denote the set of $r$--tuples of Young diagrams whose total size
is no more than $n$ and congruent to $n \mod 2$.  If   $\labold \in \Gamma_n$   and $\mubold \in \Gamma_{n-1}$,  write $\mubold \leftrightarrow \labold$ if  $\mubold$ is obtained from $\labold$ by either adding one box to, or removing one box from, one component of $\labold$.     For $\labold \in \Gamma_n$,  an {\em up--down tableau  of length $n$ and shape $\labold$} is a sequence
$
(\labold\spp 0,  \labold\spp 1,  \dots, \labold\spp n),
$
where  
\begin{enumerate}
\item
for all $i$,   $\labold\spp i \in \Gamma_i$,  
\item
$\labold\spp0 = (\emptyset, \dots, \emptyset)$,  
\item $\labold\spp n = \labold$,  and 
\item for each $i$,
$\labold\spp i \leftrightarrow \labold\spp {i+1}$.
\end{enumerate}

Let $\mathcal T(n)$ denote the set of up--down tableaux of length $n$,  and
$\mathcal T(n, \labold)$ the set of up--down tableaux of length $n$ and shape $\labold$.

\begin{definition}
Let $\labold$ be an $r$--tuple of Young diagrams.  Label the nodes (or cells, or boxes) of $\labold$ by triples $(j, x, y)$, where $j$ is the component ($1 \le j \le r$),  $x$ is the row coordinate, and $y$ is the column coordinate of the node.  Define the {\em multiplicative content} of a node $\alpha = (j, x, y)$ by $\tilde c(\alpha) = u_j q^{2 (y - x)}$.   

An {\em addable node} of $\labold$ is a node which can be added to one component to give a new
$r$--tuple of Young diagrams.  A {\em removable node} is a node which can be removed with the result still being an $r$--tuple of Young diagrams. 

 For an addable node of 
$\labold$,  set $b(\alpha, \labold) = \tilde c(\alpha)$,  and for  removable node, set
$b(\alpha, \labold) = \tilde c(\alpha)\inv$.

If $T = (\labold\spp 0,  \labold\spp 1,  \dots, \labold\spp n)$ is an up--down tableau of length $n$ and $j \le n$,  let $b(j, T) =  b(\alpha_j, \labold\spp {j-1})$,  where $\alpha_j$ is the node that is added or removed in passing from  $\labold\spp {j-1}$  to $ \labold\spp {j}$.
\end{definition}

The total number of addable and removable nodes of a single Young diagram is always odd, so the the total number of addable and removable nodes of an $r$--tuple of Young diagrams has the same parity as $r$.

In the following we work over the field of fractions  $F$ of the universal admissible integral domain $\overline R_+$.  Since $\overline R_+$ imbeds in the universal $u$-admissible ring, according to Theorem \ref{theorem:  universal admissible rings}, 
we will now denote the parameters of  $\overline R_+$ by $\rhobold$, $\qbold$,
$\deltabold_j$,  and $\ubold_i$. Recall that $\qbold$,  $\ubold_1, \dots, \ubold_r$  are algebraically independent over $\Z$.

\begin{theorem} \label{theorem:  generic structure}

 Let $F$ denote the field of fractions of the universal admissible integral domain  $\overline R_+$.
Write $\bmw{n}$  for $\bmw{n, F, r}(\ubold_1,  \dots, \ubold_r)$,  
$\kt n$ for $\kt{n, F, r}(\ubold_1,  \dots, \ubold_r)$,  
and $\hec{n}$  for $\hec{n, F, r}(\qbold^2; \ubold_1,  \dots, \ubold_r)$.  The following statements hold for all $n \ge 0$.
\begin{enumerate}
\item $\varphi: \bmw {n} \rightarrow \kt {n}$ is an isomorphism.
\item The Markov trace $\eps$ on $\kt {n}$ is non--degenerate.
\item   $\bmw{n}$ is split semisimple, and $\bmw{n} \cong I_n \oplus \hec{n}$.
\item The minimal ideals  of $\bmw{n}$  are labeled by $\Gamma_n$; more precisely,  the minimal ideals of $I_n$  are labeled by  $\{\labold \in \Gamma_n :|\labold| < n\}$ and the minimal ideals in $W_n/I_n$ are labeled by $\{\labold \in \Gamma_n :|\labold| = n\}$.
\item  The branching diagram for $W_{n-1} \subseteq W_n$  has a single edge connecting $\mubold \in \Gamma_{n-1}$  and $\labold \in \Gamma_{n}$ if  $\mubold \leftrightarrow \labold$, and no edges otherwise.
\item  The set of paths on the branching diagram for the sequence
$W_0 \subseteq W_1 \subseteq \cdots \subseteq W_n$  is the set of up--down tableaux of length $n$.  
\item Let $p_T$  be the path idempotent in $W_n$ corresponding to an up--down tableau $T$ of length $n$.  Then
 $p_T$ is a common eigenvector for $y_1, y_2, \dots, y_n$,   $$y_k \,p_T = 
 p_T \,y_k =  b(k, T) \,p_T.$$
\end{enumerate}
\end{theorem}

\begin{proof}  The assertions of the theorem hold trivially for $n = 0$, and      
the assertions for $n = 1$ are implied  by Lemma 
\ref{lemma: isomorphism and semisimplicity for one strand}. (We label the one minimal ideal of $\bmw {0, F} \cong F$  by the  $r$--tuple of empty diagrams,
$\bm \emptyset = (\emptyset, \dots, \emptyset)$.   We label the minimal ideal $F p_j$ of $\bmw {1}$  by the the $r$--tuple $\square_j$,  with the  Young diagram $\square$ 
 in the $j$--th position,  and the empty diagram $\emptyset$  in all other positions.
 Assertion (2) holds for $n = 1$ since $\eps(p_j) = \gamma_j/\delta_0$ by Lemma \ref{lemma: isomorphism and semisimplicity for one strand} and $\gamma_j \ne 0$ by 
Remark \ref{remark: new formulas for gammas}.
 Assertion (7) holds for $n = 1$ because $y _ 1 p_j = u_j p_j$,  and $u_j$ is the multiplicative content of the node added to $(\emptyset, \dots, \emptyset)$ to get
 $\square_j$.
 )

We proceed by induction on $n$.    Fix $n \ge 1$, and  suppose that the assertions hold for $\bmw f$  for $f \le n$.

We have a  conditional expectation $\eps_n : \bmw{n} \rightarrow \bmw{n-1}$ that preserves the non-degenerate trace on $\bmw{n}$ and $\bmw{n-1}$;  that is,   $\eps\circ\eps_n = \eps$.  (The conditional expectation results from the isomorphism $\bmw f \cong \kt f$ for $f \le n$ and the conditional expectation on the cyclotomic  Kauffman tangle algebras; see Section \ref{subsection: Inclusion, conditional expectation, and trace}.)
Moreover,  the conditional expectation is implemented by the idempotent $e =  (1/\deltabold_0) e_n \in \bmw{n+1}$;   that is  $e b e =  \eps_n(b) e$  for $b \in \bmw{n}$.  Therefore,  by Wenzl's Theorem \ref{Theorem: Wenzl extensions} and the induction hypothesis, 
$W_n e_n W_n \subseteq W_{n+1}$  is split semisimple,  with minimal ideals labeled by $\Gamma_{n-1}$.
By Lemma \ref{lemma:  I n plus 1 lemma},   $W_n e_n W_n = I_{n+1}$.    

Since $I_{n+1} = W_n e_n W_n$ has an identity element $z$,  necessarily a central idempotent in $W_{n+1}$,  it follows that $I_{n+1}$ has 
a complementary ideal $J_{n+1}$;   moreover,  $J_{n+1}  \cong \hec{n+1}$, according to Lemma  \ref{lemma:  hecke quotient of cyclotomic bmw}.  By  
Theorem  \ref{theorem:  cyclotomic hecke split semisimple} and
the algebraic independence of 
$\qbold, \ubold_1, \dots, \ubold_r$ over $\Z$,   $H_{n+1}$  is split semisimple.  Thus  
$\bmw{n+1}$ is split semisimple and $\bmw{n+1} \cong I_{n+1} \oplus \hec{n+1}$.  This proves statement (3)  for
$\bmw{n+1}$.

To prove statements (4),  (5) and (6)  for $\bmw{n+1}$,  we have to determine the branching diagram for $\bmw{n} \subseteq \bmw{n+1}$.    We know that the multiplicity matrix for the map   $x \mapsto  x z$,  from $\bmw{n}$  to $I_{n+1}$, 
is the transpose of the inclusion matrix for $\bmw{n-1} \subseteq \bmw{n}$, by  Theorem \ref{Theorem: Wenzl extensions}.    The map   $x \mapsto x(1 -z)$, from $I_n $  to  $J_{n+1}$,    is zero, since
$I_n \subseteq I_{n+1}$.    The multiplicity matrix for the map  $x \mapsto x(1 -z)$, from $J_n $  to  $J_{n+1}$
is the same as that for the inclusion $\hec{n}  \subseteq \hec{n+1}$,  which is given by Theorem
\ref{theorem:  cyclotomic hecke split semisimple}.    Statements (4),  (5), and (6)  follow from these observations.

It remains to prove statements (1), (2)  and (7)  for $\bmw {n+1}$  and $\kt {n+1}$.

Recall that for any up--down tableau $T$  of length $k \le n+1$,  
$T'$  denotes the truncation of $T$ of length
$k-1$.  If $T$ is an up--down tableau of length $k$ and shape $\mubold$, and if
$\labold \in \Gamma_{k+1}$  satisfies $\mubold \leftrightarrow \labold$,  then
$T + \labold$  denotes the extension of $T$ of length $k+1$ and shape $\labold$.
Let $p_T$ denote the path idempotent in $\bmw {k}$ corresponding to an up--down tableau $T$ of length $k$, and let 
$P_T = \varphi(p_T)$  be its image in $\kt k$.

\begin{lemma} \label{lemma:  action of y's on path idempotents}
Assume that assertions (1), (2), and  (7) hold  for $\bmw k$ and $\kt k$ for all  $k \le n$.
If  $T$ is an up--down tableau of length $n+1$, then  for all $k \le n+1$,
 $$y_{k}\, p_T = 
 p_T\, y_{k}  = b(k, T)\, p_T.$$
 That is,  assertion (7)  also holds for $\bmw {n+1}$.
\end{lemma}

\noindent
{\em Proof of Lemma \ref{lemma:  action of y's on path idempotents}.}   
Let $T$ be an up--down tableau of length $n+1$.
If $k \le n$,  then
$$
y_k p_T = y_k p_{T'} p_T =  b(k, T') p_{T'} p_T = b(k, T) p_T,
$$
 and similarly for $p_T y_k$.
 
 It remains to show that
 $
 y_{n+1} p_T = p_T y_{n+1} =  b(n+1, T)  p_T.
 $
 
 Write $T = (\bm \emptyset, \labold \spp 1, \dots, \labold\spp {n-1}, \mubold, \labold)$.  
 In case $|\labold| = n+1$,  $p_T$ belongs to the ideal $J_{n+1}$, which is isomorphic to the cyclotomic Hecke algebra $H_{n+1}$.  The result then follows from formulas for the irreducible representations of the cyclotomic Hecke algebras, see ~\cite{ariki-koike, ariki-book}.
 
 We assume now that $|\labold| \le n-1$.
 Consider the minimal ideal $W_{n, \mubold}$ of $W_n$ corresponding to $\mubold$.  This ideal is isomorphic to a full matrix algebra over $F$, and the path idempotents $p_U$
 ($U \in \mathcal T(n, \mubold)$)  are a complete family of mutually orthogonal minimal idempotents in $\bmw {n, \mubold}$.  Hence
 $\bmw {n, \mubold}$ has a system of matrix units $e_{S, U}$  ($S, U \in \mathcal T(n, \mubold)$)  such that 
 $
 p_U =  e_{U, S}\, p_S \,e_{S, U}.
 $
 Thus for any $S \in \mathcal T(n, \mubold)$,  we can write
 $$
 \begin{aligned}
 p_T &= p_{T'} z\spp {n+1}_\labold = e_{T', S}\, p_S\, e_{S, T'}  z\spp {n+1}_\labold \\
 &= e_{T', S}\, p_S  z\spp {n+1}_\labold \, e_{S, T'} =  e_{T', S}\, p_{S + \labold} \, e_{S, T'}.
 \end{aligned}
 $$
Suppose that there exists some $S \in \mathcal T(n, \mubold)$ such that
 $$
  y_{n+1} p_{S + \labold} = p_{S + \labold} y_{n+1} =  b(\alpha, \mubold)  p_{S + \labold},
 $$
 where $\alpha$ is the node added or removed in passing from $\mubold$ to $\labold$.
Then we have
$$
 \begin{aligned}
 y_{n+1} p_T &=  y_{n+1} e_{T', S}\, p_{S + \labold} \, e_{S, T'}
 = e_{T', S}\, y_{n+1} p_{S + \labold} \, e_{S, T'} \\ &=
 b(\alpha, \mubold)   e_{T', S}\, p_{S + \labold} \, e_{S, T'}
 =b(\alpha, \mubold) p_T,
  \end{aligned}
$$
and similarly for $ p_T y_{n+1}$.  Therefore, it suffices to prove the result when $T$ 
has the form
 $T = (\bm \emptyset, \labold \spp 1, \dots, \labold, \mubold, \labold).$ 
Let $T'$ be the truncation of $T$ of length $n$ and $T''$ the truncation of length $n-1$.
Then $p_{T''}$ is a minimal idempotent in the minimal ideal $\bmw {n-1, \labold}$ of
$\bmw {n-1}$.
Let $e = (1/\delta_0) e_n$, so that $e b e = \eps_n(b) e$ for $b \in \bmw n$.  Then by
Remark \ref{remark:  on Wenzl extensions},   $p_{T''} e$ is a minimal idempotent in $\bmw {n+1, \labold}$.

Consider the element $p_{T'} \,e \, p_{T'}$.  We have
$$
\begin{aligned}
(p_{T'} \,e \, p_{T'})^2 &= p_{T'} \,e  \eps_n(p_{T'})  \, p_{T'}   = (\eps(p_{T'})/\eps(p_{T''})) p_{T'} \,e \, p_{T''}  p_{T'}  \\&= (\eps(p_{T'})/\eps(p_{T''})) p_{T'} \,e \,  p_{T'},
\end{aligned}
$$
using Lemma \ref{lemma:  conditional expectation of path idempotents}, so $p_{T'} \,e \, p_{T'}$ is an essential idempotent.  Write $\kappa =
 \eps(p_{T'})/\eps(p_{T''})$ and $f = \kappa\inv p_{T'} \,e \, p_{T'}$.   (We are using the assumption that the trace is non--degenerate on on $\bmw n$ and $\bmw {n-1}$, which follows from assertions (1) and (2) for these algebras.)

 We have
 $f = \kappa\inv  p_{T'}( p_{T''}\,e) \, p_{T'}$,  which shows that $f$ is in
 $\bmw {n+1, \labold}$.    But then  $f = z \spp {n+1}_\labold f  z \spp {n+1}_\labold =
 ( z \spp {n+1}_\labold p_{T'}) \,f   \,(p_{T'} z \spp {n+1}_\labold) = p_T \,f \,p_T$.
 Therefore $f$ is a multiple of, and hence equal to, the minimal idempotent $p_T$.
Now we can compute:
$$
\begin{aligned}
y_{n+1} p_T &= \kappa\inv y_{n+1} p_{T'} \,e\, p_{T'} = \kappa\inv  p_{T'} y_{n+1}\,e\, p_{T'}  \\
&= \kappa\inv  p_{T'} y_{n}\inv \,e\, p_{T'} =  \kappa\inv y_{n}\inv  p_{T'}  \,e\, p_{T'} \\
&= y_n\inv p_T = b(\alpha, \mubold) p_T,
\end{aligned}
$$
and similarly for $p_T y_{n+1}$.
This completes the proof of Lemma \ref{lemma:  action of y's on path idempotents}.

Recall the quantities $Z_{k+1}(t)$  and  $Q_{k+1}(t)$ from Section \ref{subsection:  more generating functions}.  They are  formal power series with coefficients in the center of $\kt {k}$.  
 We want to find an expression for $Q_{k+1}(t) p_T$,   where $T$ is an up--down tableau of length $k$.  
 
 Let $\labold$ be any $r$--tuple of Young diagrams.     Define
 \begin{equation}  \label{equation:  expression for tilde Q of t and lambda}
 \tilde Q(t, \labold) = p\  A(t) \prod_\alpha  \frac{t - b(\alpha, \labold)\inv}{t - b(\alpha, \labold)},
 \end{equation}
where $p = \prod_j u_j$,   $A(t)$ is as in  (\ref{equation:  formula for Z1 in u admissible case}),  and $\alpha$ runs over all addable and removable nodes of $\labold$.  Define
\begin{equation}  \label{equation: definition of tilde Z of t and lambda}
\tilde Z(t, \labold) =  \tilde Q(t, \labold) + \rho\inv/(q\inv - q)  +  t^2/(t^2 -1).
\end{equation}

Note that if $\bm \emptyset$ is the $r$--tuple of empty Young diagrams,  then
 \begin{equation}  \label{equation:  expression for tilde Q of t and empty}
 \tilde Q(t, \bm \emptyset) = p\  A(t) \prod_j  \frac{t - u_j\inv}{t - u_j} = A(t) \,G(t\inv) = Q_1(t), 
 \end{equation}
 using  (\ref{equation: expression for Q1 of t}) and thus $\tilde Z(t, \bm \emptyset) = Z_1(t)$.

\begin{lemma} \label{lemma:  recursion for tilde Q of t and lambda}
Suppose $\mubold$ and $\labold$ are $r$--tuples of Young diagrams and $\labold$ is obtained from $\mubold$ by adding a single node $\alpha$.   Let $\tilde c = \tilde c(\alpha)$ be the multiplicative content of this node.  Then
$$
\tilde Q(t, \labold) = \tilde Q(t, \mubold) \frac{(t - \tilde c)^2 (t - q^{-2} \tilde c\inv)(t - q^{2} \tilde c\inv)}
{(t - \tilde c\inv)^2 (t - q^{-2} \tilde c)(t - q^{2} \tilde c)}
$$
\end{lemma}

\noindent
{\em Proof of Lemma \ref{lemma:  recursion for tilde Q of t and lambda}.}   Suppose $\alpha = (j, x, y)$,  so $y-1 =\mu_x\spp j$,  the length of the the $x$--th row of the $j$--th component of $\mubold$.   There are several cases to consider according to the relative size of
$\mu_{x-1}\spp j$,  $\mu_{x}\spp j$, and $\mu_{x+1}\spp j$.  For example,  suppose
$\mu_{x}\spp j = \mu_{x+1}\spp j$,  and  $\mu_{x}\spp j + 1 <  \mu_{x-1}\spp j$.  Then adding the node  $\alpha$:
\begin{enumerate}
\item  eliminates the addable node at $\alpha$ and produces a removable node at $\alpha$, and
\item  produces two new addable nodes at $(j, x, y+1)$  and $(j, x+1, y)$.
\end{enumerate}
The node at $\alpha$ contributes $({(t- \tilde c)}/{(t - \tilde c\inv)})^2$ to the ratio
${\tilde Q(t, \labold)}/{\tilde Q(t, \mubold)}$.  The two new addable nodes contribute
$$\displaystyle \left (\frac{t- q^{-2} \tilde c\inv}{t - q^2 \tilde c}\right ) \left (\frac{t- q^{2} \tilde c\inv}{t - q^{-2} \tilde c}\right).
$$  
Other cases can be analyzed similarly.  This completes the proof of Lemma  \ref{lemma:  recursion for tilde Q of t and lambda}.

\begin{lemma} \label{lemma: action of Qk of t on path idempotent}
 Let $T$ be an up--down tableau of length $k \le n+1$ and shape $\labold$  Then
$$
Q_{k+1}(t) P_T  = \tilde Q(t, \labold) P_T.
$$
and
$$
Z_{k+1}(t) P_T  = \tilde Z(t, \labold) P_T.
$$
\end{lemma}

 \noindent
{\em Proof of Lemma  \ref{lemma: action of Qk of t on path idempotent}.}  We do this by induction on $k$.  There is a unique tableau of length $0$,  and the corresponding idempotent is the identity.   So the statement for $k = 0$  reads $Q_1(t) = \tilde Q(t, \bm \emptyset)$, which was observed in 
 (\ref{equation:  expression for tilde Q of t and empty}).  Let $T$ be an up--down tableau of length $k \ge 1$  and suppose the assertion holds for all shorter up--down tableaux.  
Let $\labold$ denote the shape of $T$  and $\mubold$ the shape of $T'$.
Then
$$
\begin{aligned}
Q_{k+1}(t) P_T &=   Q_k(t) \left (  \frac{(t-Y_k)^2 (t - q^{-2} Y_k\inv) (t- q^2 Y_k\inv)}
{(t-Y_k\inv)^2 (t - q^{-2} Y_k) (t- q^2 Y_k)} \right )  P_{T'}  P_T \\
&=  Q_k(t)  P_{T'} \left (  \frac{(t-Y_k)^2 (t - q^{-2} Y_k\inv) (t- q^2 Y_k\inv)}
{(t-Y_k\inv)^2 (t - q^{-2} Y_k) (t- q^2 Y_k)} \right )   P_T.
\end{aligned}
$$
By the induction hypothesis and Lemma \ref{lemma:  action of y's on path idempotents}, this is equal to 
$$
\begin{aligned}
\tilde Q(t, \mubold) &P_{T'}  \left (  \frac{(t-b)^2 (t - q^{-2} b\inv) (t- q^2 b\inv)}
{(t-b\inv)^2 (t - q^{-2} b) (t- q^2 b)} \right )   P_T \\
& = \tilde Q(t, \mubold)   \left (  \frac{(t-b)^2 (t - q^{-2} b\inv) (t- q^2 b\inv)}
{(t-b\inv)^2 (t - q^{-2} b) (t- q^2 b)} \right )   P_T, 
\end{aligned}
$$
where $b = b(\alpha, \mubold)$.  By Lemma \ref {lemma:  recursion for tilde Q of t and lambda}, this equals $\tilde Q(t, \labold) P_T$.   (In case $\labold \subseteq \mubold$,  reverse the roles of 
$\labold$ and $\mubold$ in the Lemma.)    The conclusions of Lemma \ref{lemma: action of Qk of t on path idempotent} follow from this.

Let $U$ be an up--down tableau of length $n$.
Note that $P_{U}$ is a minimal idempotent in $\kt {n}$ and $\eps(P_{U}) \ne 0$  by claims (1) and (2)  of the induction hypothesis.   We have $P_U = \sum_T  P_T$,  where  $T$  ranges over all
up--down tableaux of length $n+1$  such that $T' = U$.  It follows that there exists some
$T$ such that $\eps(P_T) \ne 0$.   Fix such a $T$.

 Recall from Lemma \ref{lemma:  conditional expectation of path idempotents} that  for any up--down tableau $S$  of length $n+1$  such that $S' = U$,  we have
 \begin{equation} \label{equation:  formula for cond expec of path idempotent2}
 \eps_{n+1}(P_S) =   (\eps(P_S)/ \eps(P_{U})) P_{U}.
 \end{equation}
 Next,  for any such  up--down tableau $S$,  we claim
\begin{equation}  \label{equation:  formula involving PT En PS En PT}
\eps(P_T E_{n+1} P_S E_{n+1} P_T) =  \frac{\eps(P_T) \eps(P_S)}{\eps(P_{U})}.
\end{equation}
In fact,
$$
\begin{aligned}
\eps(& P_T E_{n+1} P_S E_{n+1} P_T)  = \eps( E_{n+1} P_S E_{n+1} P_T)   \\
&= \deltabold_0\  \eps( E_{n+1}  \eps_{n+1}(P_S)  P_T)  
= \deltabold_0 \frac{\eps(P_S)}{\eps(P_{U})} \eps (E_{n+1} P_{T'}  P_T) \\
&= \deltabold_0 \frac{\eps(P_S)}{\eps(P_{U})}  \eps (E_{n+1}  P_T)  
=  \frac{\eps(P_T) \eps(P_S)}{\eps(P_{U})}, 
\end{aligned}
$$
where we have used Lemma \ref {lemma:  formulas involving cond exp and En}
and  (\ref{equation:  formula for cond expec of path idempotent2}).

Now consider the formal power series:
\begin{equation}  \label{calculation of weights 1}
\begin{aligned}
E_{n+1} (t - Y_{n+1})\inv E_{n+1}&= E_{n+1} (1/t) \sum_{a \ge 0} Y_{n+1}^a t^{-a}  E_{n+1} \\
&= (1/t) Z_{n+1}(t) E_{n+1}
\end{aligned}
\end{equation}
 Let $U$ and $T$ be as above,  so $T' = U$ and $\eps(P_T) \ne 0$.  Denote the shape of $U$ by
 $\mubold$.  
Multiply both sides of  \ref{calculation of weights 1} by $P_T$,  and use Lemma \ref{lemma: action of Qk of t on path idempotent} to get
\begin{equation}  \label{calculation of weights 2}
\begin{aligned}
P_T E_{n+1} (t - Y_{n+1})\inv E_{n+1}  P_T= (1/t) \tilde Z(t, \mubold)  P_T E_{n+1}  P_T.
\end{aligned}
\end{equation}
Now take the trace on both sides,
\begin{equation}  \label{calculation of weights 3}
\begin{aligned}
\eps(P_T E_{n+1} (t - Y_{n+1})\inv E_{n+1}  P_T)= (1/t) \tilde Z(t,\mubold)  \delta_0\inv \eps(P_T).
\end{aligned}
\end{equation}
Now we focus on the left side of the equation.  Note that for $S$ an up--down tableau of length
$n+1$,  $P_S E_{n+1} P_T =  P_S P_{S'} P_{T'}  E_{n+1} P_T$,  and this is non--zero only if 
$S' = T' = U$.  Therefore,  the left side equals
$$
\sum_S \eps(P_T E_{n+1} (t - Y_{n+1})\inv P_S E_{n+1}  P_T),
$$
where the sum is over all $S$ such that $S' = U$.  By Lemma \ref{lemma:  action of y's on path idempotents},  and  (\ref{equation:  formula involving PT En PS En PT})
 this equals
\begin{equation} \label{calculation of weights 4}
\begin{aligned}
\sum_S  &(t - b(n+1, S))\inv \eps(P_T E_{n+1} P_S E_{n+1}  P_T) \\&= \sum_S  (t - b(n+1, S))\inv \frac{\eps(P_S) \eps(P_T)}{\eps(P_U)} \\
&= \sum_\alpha (t - b(\alpha, \mubold))\inv \frac{\eps(P_{S(\alpha)}) \eps(P_T)}{\eps(P_U)},
\end{aligned}
\end{equation}
where the last sum is over all addable and removable nodes of $\mubold$, and
$S(\alpha)$ is the extension of $U$ by $\alpha$;  that is,  if 
$
U = (\lambdabold \spp 0, \dots, \mubold),
$
then
$
S(\alpha) = (\lambdabold \spp 0, \dots, \mubold, \labold(\alpha)),
$
where $\labold(\alpha)$ is obtained from $\mubold$ by addition or removal of the node $\alpha$.
  Comparing 
\ref{calculation of weights 3} and \ref{calculation of weights 4} and cancelling $\eps(P_T)$, we have
\begin{equation} \label{calculation of weights 5}
\begin{aligned}
  (1/t) \tilde Z(t,\mubold) =  \deltabold_0 \sum_\alpha (t - b(\alpha, \mubold))\inv \frac{\eps(P_{S(\alpha)})}{\eps(P_U)},
\end{aligned}
\end{equation}
Since $\qbold$, $\ubold_1$,  \dots, $\ubold_r$  are algebraically independent over
$\Z$,  if $\alpha$ and $\beta$ are two different addable or removable nodes of $\mubold$, then $\tilde c(\alpha) \ne \tilde c(\beta)\powerpm$; hence $b(\alpha, \mubold) \ne b(\beta, \mubold)$.  
Let $S$ be any particular extension of $U$,   $S = S(\beta)$,  where $\beta$ is an addable or removable node of $\mubold$;  let $b =  b(\beta, \mubold)$.   Take the residue at $t = b$  on both sides of  \ref{calculation of weights 5}.  On the left side, we get
\begin{equation} \label{calculation of weights 6}
\begin{aligned}
  \Res (1/t) \tilde Z(t,\mubold)_{| t = b} &=   \Res (1/t) \tilde Q(t,\mubold)_{| t = b} \\
  &= p b\inv A(b) (b-b\inv)\prod_{\alpha \ne \beta} \frac{b - b(\alpha, \mubold)\inv}{b - b(\alpha, \mubold)}, \\
\end{aligned}
\end{equation}
while on the right side we get 
  $\deltabold_0  {\eps(P_S)}/{\eps(P_U)}$.  Therefore,
  \begin{equation} \label{calculation of weights 7}
\begin{aligned}
\eps(P_S) =   \deltabold_0\inv \,p\, b\inv A(b)(b-b\inv) \prod_{\alpha \ne \beta} \frac{b - b(\alpha, \mubold)\inv}{b - b(\alpha, \mubold)} \eps(P_U).
\end{aligned}
\end{equation}
If $r$ is odd,  $A(b) (b - b\inv) = \displaystyle \frac{b-b\inv}{q - q\inv} + 1$, and $p = \rhobold$; see  (\ref{equation:  formula for Z1 in u admissible case}).  Therefore,
  \begin{equation} \label{calculation of weights 8}
\begin{aligned}
\eps(P_S) =  \deltabold_0\inv \,\rhobold\, b\inv (\frac{b-b\inv}{q - q\inv} + 1) \prod_{\alpha \ne \beta} \frac{b - b(\alpha, \mubold)\inv}{b - b(\alpha, \mubold)} \eps(P_U).
\end{aligned}
\end{equation}
If $r$ is even,  $A(b) (b - b\inv) = \displaystyle \frac{q\inv b - q b\inv}{q - q\inv}$, and
$p = q\inv \rhobold$.  Therefore,
  \begin{equation} \label{calculation of weights 9}
\begin{aligned}
\eps(P_S) =  \deltabold_0\inv \,\rhobold \,q\inv b\inv (\frac{q\inv b-q b\inv}{q - q\inv} ) \prod_{\alpha \ne \beta} \frac{b - b(\alpha, \mubold)\inv}{b - b(\alpha, \mubold)} \eps(P_U).
\end{aligned}
\end{equation}
The algebraic independence of  $\qbold$, $\ubold_1$,  \dots, $\ubold_r$ over $\Z$ implies that $\eps(P_S) \ne 0$.

The computation of $\eps(P_S)$ given here was  adapted from   ~\cite{blanchet-beliakova}, Theorem 7.

We have shown that for every up--down tableau $S$ of length $n+1$,  $\eps(P_S) = \eps(\varphi(p_S)) \ne 0$.   Hence $\eps$ is a non--degenerate trace on $\bmw {n+1}$.  Moreover,  this shows that $\varphi$ is not zero on any minimal ideal of $\bmw {n+1}$,  and hence $\varphi$ is injective on $\bmw {n+1}$.   This proves points (1) and (2)  for
$\bmw {n+1}$,  and completes the proof of the theorem.
\end{proof}

\begin{remark} (A semisimplicity criterion)  Let $K$ be a field with $u$--admissible parameters.  The inductive proof of the theorem applies to $\bmw {n, K, r}$  as long as the weights of the Markov trace remain non--zero.  For this, what is needed is that for all $r$--tuples of Young diagrams $\mubold $ of size $\le n$,  and for all addable or removable nodes $\alpha$ and $\beta$ of $\mubold$,  with $b = b(\beta, \mubold)$ and $b' = b(\alpha, \mubold)$, we have
\begin{equation*}
A(b)(b - b\inv) \ne 0, \ \text{and}\  b \ne (b')^{\pm 1}.
\end{equation*}
For this, it suffices that for all $j$,  $u_j \ne \pm q^{\pm (2k -1)}$, and for all $j, j'$, 
$u_j/u_{j'} $ and $u_j u_{j'}$  are not equal to $q^{2k}$  for $k \le n$.  In particular, we have the following result:
\end{remark}

\begin{corollary}  Let $K$ be a field with $u$--admissible parameters.  Suppose that for all 
$j, j'$,  none of the quantities $\pm u_j$,  $u_j/u_{j'}$  and $u_j u_{j'}$  is equal to an integer power of $q$.  Then $\bmw{n, K, r}$ is split semisimple for all $n$.
\end{corollary}

\section{The $\Z_r$--Brauer algebra and the dimension of $\bmw{n, F,r}$.}
Recall that for $\labold \in \Gamma_n$,  the set of up--down tableaux of length $n$ and shape $\labold$ is denoted by  $\mathcal T(n, \labold)$.   The dimension of 
the simple $\bmw {n, F, r}$ module labelled by $\labold \in \Gamma_n$ is 
$|\mathcal T(n, \labold)|$, and 
the dimension of $\bmw {n, F,r}$ over $F$ is the $\sum_{\labold \in \Gamma_n} |\mathcal T(n, \labold)|^2$.  It should possible to compute this quantity by an appropriate Shensted type correspondence,  but it is more convenient to exhibit another sequence of split semisimple algebras  $D_0 \subseteq D_1 \subseteq \cdots \subseteq D_n$
(over some field)  whose dimensions are known and which have the same branching diagram.   

For this purpose, we introduce the $\Z_r$--Brauer algebras.
The $\Z_r$--  and $\Z$--Brauer algebras
 were considered in  ~\cite{H-O2},  ~\cite{Rui-Yu}, and   ~\cite{GH1}. 
 The $G$--Brauer algebras for an arbitrary
 abelian group $G$  were studied in ~\cite{G-Brauer}.

We recall that a Brauer diagram is a tangle diagram in the plane,  in which information about over-- and under--crossings is ignored:
An $(n,n)$--{\em Brauer diagram}  is a figure in the rectangle  $R = I \times I$ 
consisting
of
\begin{enumerate}
\item  $n$ distinguished points $\{\p 1,  \dots, \p n\}$  on the upper boundary 
$I \times \{1\}$ of $R$  and $n$ distinguished points   $\{ \pbar 1,
\dots \pbar n\}$ on the lower boundary $I \times \{0\}$ of $R$.
\item $n$  curves connecting  the points $\{\p 1,  \dots, \p n, \pbar 1,
\dots \pbar n\}$ in pairs, with each curve intersecting the boundary of $R$ only at its two endpoints.
\end{enumerate}
A $\Z_r$--Brauer diagram is
a Brauer diagram in which each curve  (strand) is endowed with an orientation and labeled by an element of the group $\Z_r$.   Two labelings are regarded as the same if the orientation of a strand is reversed and the group element associated to the strand is inverted.  Note that the number of $\Z_r$--Brauer diagrams is
$
 r^n (2n-1)!!.
$

Let $A$ be an integral domain containing  elements $\vt_j$, for $0 \le j \le \lfloor r \rfloor$  with
 $\vt_0$  invertible.
 The  $\Z_r$--Brauer algebra $\br{n,  A, r}$ is defined as follows.
As an $A$--module,  $\br{n, A, r}$ is the free  $A$--module
 with basis the set
of  $\Z_r$--Brauer diagrams with $n$ strands.      The product of two  $\Z_r$--Brauer diagrams
is defined to be a certain multiple of another  $\Z_r$--Brauer
diagram, determined as follows.

Given two $\Z_r$--Brauer diagrams $a, b$, first ``stack" $b$ over $a$ (as
for tangle diagrams and ordinary Brauer diagrams).   In the resulting
``tangle'',  the strands    are composed of one or more
strands from $a$  and $b$.   Give each composite strand $s$  an orientation,  arbitrarily.
Make the orientations of the components of $s$  from $a$  and $b$  agree with the
orientation of $s$, changing the signs of the $\Z_r$--valued labels accordingly;  the label
of  $s$  is then the sum of the labels of its components from $a$ and $b$.
 For each $j$ ($0 \le j \le \lfloor r \rfloor$), let  $m_j$  be the number of  closed loops with label  equal to $\pm[j]$.
Let $c$ be the 
$\Z_r$--Brauer diagram obtained by removing all the closed curves.  Then
$
a b =  (\prod_{j}   \vt_{j}^{m_{j}} )\  c .
$
Extend the multiplication to a bilinear product on $\br{n, A, r}$.
One can easily check that the multiplication is associative.  
Note that
the subalgebra generated by $\Z_r$--Brauer diagrams with only vertical
strands is isomorphic to the group algebra of the wreath product $\Z_r \wr \S_n$, so the $\Z_r$--Brauer algebra can be considered as a sort of wreath product of $\Z_r$ with the ordinary Brauer algebra.

One has inclusion maps $\iota: \br{n-1, A, r} \rightarrow \br{n, A, r}$  defined on the level of $\Z_r$--Brauer diagrams  by
adding a vertical strand on the right labeled by $0$.
Moreover, one has a conditional expectation $\eps_n : \br{n, A, r} \rightarrow
\br{n-1, A, r}$ defined as follows:  First define a map $\cl_n$,  from
$\Z_r$--Brauer diagrams with $n$--strands to $\Z_r$--Brauer diagrams with $(n-1)$--strands, 
 by joining the rightmost
pair of vertices  $\pbar n, \p n$    of a diagram $d$  by a
new strand, with label $0$:
$$
{\rm cl}_n: \quad \inlinegraphic{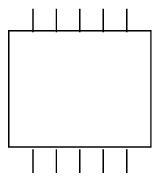} \quad \mapsto \quad 
\inlinegraphic{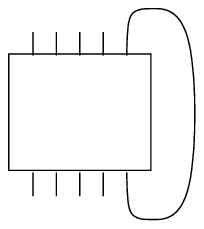}
$$
Determine the $\Z_r$--valued label of any concatenated strand in the resulting
diagram by the same rule as in the multiplication of $\Z_r$--Brauer diagrams.
If the resulting diagram has a closed loop 
(which happens precisely if $d$ contains a strand connecting
$\p n$ and
$\pbar n$ with some label $\pm [k]$) then remove the closed loop and multiply
the resulting  $\Z_r$--Brauer diagram by $\vt_{k}$.  Finally,
define $\eps_n : \br{n, A, r} \rightarrow \br{n-1, A, r} $ by 
$$
\eps_n(d) =\vt_0\inv \cl_n(d).
$$
One can check that $\eps_n$ is a conditional expectation.  
Define $$\eps = \eps_1 \circ \dots \circ \eps_n : \br{n,A, r} \rightarrow
\br{0, A, r} = A.$$  Then $\eps$ is a trace.

\begin{lemma}[\cite{G-Brauer}] \label{lemma: non degenerate trace}
For 
 each $j$ ($0 \le j \le \lfloor r \rfloor$),  let $\vtbold_{j}$ be an indeterminant
 over $\C$.    
   Let 
$$K = \C(\vtbold_0^{\pm 1},  \vtbold_{1} , \dots, \vtbold_{\lfloor r \rfloor}).$$
The trace $\eps$ on $\br{n, K, r}$ is non-degenerate. 
\end{lemma}

Given this,  it is fairly straightforward to apply Wenzl's method from ~\cite{Wenzl-Brauer, Wenzl-BCD} to the $\Z_r$--Brauer algebra. (In fact,  this has been done more generally in ~\cite{G-Brauer}, in the context of $G$--Brauer algebras,  where $G$ is a finite abelian group.)  The result of this analysis is that the $\Z_r$--Brauer algebra $\br n = \br{n, K, r}$ is
split semisimple with simple modules labelled by $\Gamma_n$ (i.e., the set of $r$--tuples of Young diagrams of total size
 $\le n$ and congruent to $n \mod 2$)  and the branching diagram for the sequence  $D_0 \subseteq D_1 \subseteq \cdots \subseteq D_n$
is exactly the same as that obtained for the sequence of cyclotomic BMW algebras in 
Theorem \ref{theorem:  generic structure}.  Thus we have:

\begin{theorem} \label{theorem on dimension of cyclotomic BMW}

Let $F$ be as in Theorem \ref{theorem:  generic structure}.  The dimension of $\bmw{n, F, r}$ is \break$r^n (2n-1)!!$.
\end{theorem}

\begin{proof}  
$\dim_F(\bmw{n, F, r}) = \sum_{\labold \in \Gamma_n} |\mathcal T(n, \labold)|^2 = \dim_K(\br{n, K,r}) = r^n (2n-1)!!.
$
\end{proof}

\section{Freeness  and isomorphism with cyclotomic Kauffman tangle algebras.}  We show in ~\cite{GH2}, Proposition 3.6,   that every cyclotomic
BMW algebra \break $\bmw{n, S, r}(u_1, \dots, u_r)$  is spanned over its ground ring $S$ by
a set  $\A_r'$  of cardinality  $r^n (2n-1)!!$.
Combining this result with Theorems  \ref{theorem:  generic structure} and
 \ref{theorem on dimension of cyclotomic BMW}, we obtain the following theorem  (which has been proved independently by Wilcox and Yu.)

\begin{theorem}[Goodman and Hauschild--Mosley,  Wilcox and Yu \cite{Yu-thesis, Wilcox-Yu2}]  Let $S$ be an admissible integral domain.  Then
$\bmw{n, S, r} \cong  \kt{n, S, r}$, and $\bmw{n, S, r}$ is a free $S$--module of rank
$r^n (2n-1)!!$.
\end{theorem}

\begin{proof}  As noted above, $\kt{n S, r}$  has a spanning set $\A_r'$ of cardinality
$r^n (2n-1)!!$.  To prove the theorem, it suffices
 to show that $\varphi(\A'_r) = \B'_r$ is linearly independent in 
$\kt{n, S, r}$.  When $S = F$,  the field of fractions of the universal admissible integral domain $\overline R_+$, this follows from Theorem \ref{theorem:  generic structure} and
Theorem \ref{theorem on dimension of cyclotomic BMW}, since $\B'_r$ is a spanning set whose cardinality equals the dimension of
$\kt{n, F, r}$.   The map $x \mapsto x \otimes 1$  from
$\kt{n, \overline R_+, r}$ to $\kt{n, \overline R_+, r} \otimes_{\overline R_+}  F \cong 
\kt{n, F, r}$ is $\overline R_+$--linear and maps  $\B'_r$ to a linearly independent set in $\kt{n, F, r}$;  hence $\B'_r$ is linearly independent in
$\kt{n, \overline R_+, r}$.    Finally,  since any admissible integral domain $S$ is a quotient of  $\overline R_+$,  and $\kt{n, \overline R_+, r}$ is a free $\overline R_+$--module with basis
$\B'_r$,   it follows that
$\kt{n, S, r} \cong \kt{n, \overline R_+, r} \otimes_{\overline R_+} S$  is a free $S$--module
with basis $\B'_r$.
\end{proof}

\begin{corollary} Suppose $S$ is an admissible integral domain.
\begin{enumerate}
\item  There is a trace--preserving conditional expectation $\eps_n : \bmw{n+1, S, r} \to \bmw{n, S, r}$.
\item  The natural map $\iota: \bmw{n, S, r} \to \bmw{n+1, S, r}$ is injective.
\end{enumerate} 
\end{corollary}

\bibliographystyle{amsplain}
\bibliography{cyclotomicbmwII}

\end{document}